\documentclass[11pt]{amsart}
\usepackage{amsmath}
\usepackage{amssymb}
\usepackage{amsthm}

\ifx\shlhetal\undefinedcontrolsequence\let\shlhetal\relax\fi
\newcount\skewfactor
\def\mathunderaccent#1#2 {\let\theaccent#1\skewfactor#2
\mathpalette\putaccentunder}
\def\putaccentunder#1#2{\oalign{$#1#2$\crcr\hidewidth
\vbox to.2ex{\hbox{$#1\skew\skewfactor\theaccent{}$}\vss}\hidewidth}}
\def\name{\mathunderaccent\tilde-3 }

\newcommand{\lgx}{\ell g} 
\newcommand{\e}{{\bf e}}
\newcommand{\n}{{\bf n}}
\newcommand{\m}{{\bf m}}
\newcommand{\Lim}{{\rm Lim}}
\newcommand{\cosik}{\Xi}
\newcommand{\Gb}{\mathfrak{b}}
\newcommand{\cov}{{\rm cov}}
\newcommand{\Con}{{\rm Con}}
\newcommand{\nul}{{\rm null}}
\newcommand{\Random}{{\rm Random}}
\newcommand{\Leb}{{\rm Leb}}
\newcommand{\cI}{{\mathcal I}}

\newcommand{\cK}{{\mathcal K}}
\newcommand{\cM}{{\mathcal M}}

\newcommand{\cP}{{\mathcal P}}
\newcommand{\cT}{{\mathcal T}}
\newcommand{\bbr}{{\mathbb R}}
\newcommand{\bbR}{\bbr}
\newcommand{\bbQ}{{\mathbb Q}}
\newcommand{\true}{{\rm true}}
\newcommand{\false}{{\rm false}}
\newcommand{\tree}{{\rm tree}}
\newcommand{\Ord}{{\rm Ord}}
\newcommand{\full}{{\rm full}}

\newcommand{\vartrianglelefteq}{\trianglelefteq}
\newcommand{\XXeq}{\vartrianglelefteq}

\newcommand{\V}{{\bf V}}

\newcommand{\rest}{{\restriction}}
\newcommand{\fs}{{}^{\textstyle \omega\!>} 2}

\newcommand{\can}{{{}^{\textstyle \omega}2}}
\newcommand{\kom}{{}^{\textstyle \omega}\kappa}

\newcommand{\cf}{{\rm cf}}
\newcommand{\tka}{{}^{\textstyle \kappa}2}

\newcommand{\dom}{{\rm dom}}

\newcommand{\forces}{\Vdash} 
\newcommand{\lesdot}{\mathrel{\mathord{<}\!\!\raise 0.8
pt\hbox{$\scriptstyle\circ$}}}  

\newcommand{\lesdoteq}{\mathrel{\mathord{\leq}\!\!\raise 0.8
pt\hbox{$\scriptstyle\circ$}}}  
\newcommand{\cH}{{\mathcal H}}
\newcommand{\cB}{{\mathcal B}}

\newcommand{\conc}{^\frown\!}

\newcommand{\supx}{{\rm supp}}

\newtheorem{theorem}{Theorem}[section] 
\newtheorem{theorem619}{Theorem}
\newtheorem{claim}[theorem]{Claim}
\newtheorem{subclaim}[theorem]{Subclaim}
\newtheorem{lemma}[theorem]{Lemma} 
\newtheorem{proposition}[theorem]{Proposition}

\newtheorem{fact}[theorem]{Fact}

\theoremstyle{definition}
 
\newtheorem{definition}[theorem]{Definition}
\newtheorem{hypothesis}[theorem]{Hypothesis}

\theoremstyle{remark}

\newtheorem{remark}[theorem]{Remark} 
\newtheorem{notation}[theorem]{Notation} 
\newtheorem{discussion}[theorem]{Discussion}

\title{Covering of the null ideal may have countable cofinality}
\author{Saharon Shelah}
\address{Institute of Mathematics\\
 The Hebrew University of Jerusalem\\
 Jerusalem 91904, Israel\\
 and  Department of Mathematics\\
 Rutgers University\\
 New Brunswick, NJ 08854, USA}
\email{shelah@math.huji.ac.il}
\urladdr{http://www.math.rutgers.edu/$\sim$shelah}
\thanks{Research supported by  ``The Israel Science Foundation'' founded
by The Israel Academy of Sciences and Humanities.  Publication no 592}

\begin{document}

\begin{abstract}
We prove that it is consistent that the covering number of the ideal of
measure zero sets has countable cofinality.
\end{abstract}

\maketitle

\section{Introduction}

In the present paper we show that it is consistent that the covering of the
null ideal has countable cofinality.  Recall that the covering number of the
null ideal (i.e. the ideal of measure zero sets) is defined as
$$
\cov(\nul)=
\min\{|{\mathcal P}|:{\mathcal P} \subseteq \nul\mbox{ and }
\bigcup_{A \in {\mathcal P}}A={\mathbb R} (= \can)\}.
$$
The question whether the cofinality of $\cov(\nul)$ is uncountable has been
raised by D.~Fremlin and has been around since the late seventies.  It
appears in the current Fremlin's list of problems, \cite{Fe94}, as problem
CO. Recall that for the ideal of meagre sets the answer is positive,
i.e. A.~Miller \cite{Mi82} proved that the cofinality of the covering of
category is uncountable. T.~Bartoszy\'nski \cite{Ba88} saw that
$\Gb<\aleph_\omega$ is necessary (see \cite[ch 5]{BaJu95} for more results
related to this problem). It should be noted that most people thought
$\cf(\cov(\nul))=\aleph_0$ is impossible

The main result of this paper is the following:

\begin{theorem}\label{A}
$\Con( \cov(\nul)=\aleph_{\omega}+{\rm MA}_{\aleph_{n}})$ for each
$n<\omega$. 
\end{theorem}

The presentation of the proof of \ref{A} sacrifices generality for hopeful
transparency. We finish by some further remarks, e.g.\ the exact cardinal
assumption for \ref{A}. We try to make the paper self contained for readers
with basic knowledge of forcing.

In a subsequent paper, \cite{Sh:619}, we deal with the question: ``can every
non-null set be partitioned to uncountably many non-null sets'',
equivalently: ``can the ideal of null sets which are subsets of a non-null
subset of ${\mathbb R}$ be $\aleph_1$-saturated''. P.~Komjath \cite{Ko}
proved that it is consistent that there is a non-meagre set $A$ such that
the ideal of meagre subsets of $A$ is $\aleph_{1}$-saturated. The question
whether a similar fact may hold for measure dates back to Ulam, see also
Prikry's thesis.  It appears as question $EL(a)$ on the Fremlin's list.  In
\cite{Sh:619} we prove the following:

\begin{theorem619}\label{B}
It is consistent that there is a non-null set $A\subseteq {\mathbb R}$ such
that the ideal of null subsets of $A$ is $\aleph_{1}$-saturated (of course,
provided that ``ZFC $+ \exists$ measurable'' is consistent).
\end{theorem619}

In \cite{Sh:619} we also prove the following.
\begin{theorem619}
\label{0.3}
It is consistent that:
\begin{enumerate}
\item[$(\oplus)$] there is a non null $A\subseteq \bbr$ such that: for
every $f: A\rightarrow \bbr$, the function $f$ as a subset of the plane
$\bbr \times \bbr$ is null
\end{enumerate}
provided that ``ZFC + there is a measurable cardinal'' is consistent.
\end{theorem619}

\begin{notation}
\label{not}
We denote: 
\begin{itemize} 
\item natural numbers by $k$, $l$, $m$, $n$ and also $i$, $j$
\item ordinals by $\alpha$, $\beta$, $\gamma$, $\delta$, $\zeta$, $\xi$
($\delta$ always limit) 
\item cardinals by $\lambda$, $\kappa$, $\chi$, $\mu$ 
\item reals by $a$, $b$ and positive real (normally
small) by $\varepsilon$.
\item subsets of $\omega$ or $^{\omega\geq}2$ or $\Ord$ by $A$, $B$, $C$,
$X$, $Y$, $Z$ but  
\item ${\mathcal B}$ is a Borel function
\item finitely additive measures by $\Xi$
\item sequences of natural numbers or ordinals by $\eta$, $\nu$, $\rho$,
\item $s$ is used for various things
\end{itemize} 
$\cT$ is as in definition \ref{T}, $t$ is a member of $\cT$.

\noindent We denote 
\begin{itemize}
\item forcing notions by $P$, $Q$, 
\item forcing conditions by $p$, $q$
\end{itemize}
and use $r$ to denote members of $\Random$ (see below) except in definition
\ref{2.1}.
\begin{itemize}
\item $\Leb$ is the Lebesgue measure (on $\{A: A\subseteq {}^\omega2\}$),
\item $\Random$ will be the family 
$$
\begin{array}{ll}
\big\{r\subseteq {}^{\omega>}2: &r\mbox{ is a subtree of }({}^{\omega >}2,
\vartriangleleft)\\
\ &\mbox{(i.e.\ non-empty subset of ${}^{\omega>}2$ closed under initial
segments)}\\  
\ &\mbox{with no }\vartriangleleft\mbox{--maximal element}\\
\ &\mbox{(so }\lim(r)=:\{\eta\in\can: (\forall n\in\omega)(\eta\rest n\in
t)\}\\
\ &\mbox{ is a closed subset of $\can$)}\\
\ &\mbox{and } \Leb(\lim(r))>0\big\}
  \end{array}$$
ordered by inverse inclusion. We may sometimes use instead
$$\{B: B \mbox{ is a Borel non-null subset of }\can\}.$$
\end{itemize}
For $\eta\in {}^{\omega >}2$, $A\subseteq {}^{\omega\geq}2$ let
$$
A^{[\eta]}=\{\nu\in A: \nu\vartrianglelefteq \eta \vee \eta \vartrianglelefteq
\nu\}.
$$
Let $\cH(\chi)$ denote the family of sets with transitive closure of
cardinality $<\chi$, and let $<^*_\chi$ denote a well ordering of
$\cH(\chi)$. 
\end{notation}

We thank Tomek Bartoszy\'nski and Mariusz Rabus for reading and
commenting and correcting.

\section{Preliminaries}
We review various facts on finitely additive measures.

\begin{definition} 
(1) ${\mathcal M}$ is the set of functions $\cosik$ from some Boolean
subalgebra $P$ of ${\mathcal P}(\omega)$ including the finite sets, to
$[0,1]_{\bbr}$ such that:
\begin{itemize}
\item $\cosik(\emptyset)=0$, $\cosik(\omega)=1$,
\item if $Y,Z \in P$ are disjoint, then $\cosik(Y \cup Z)=\cosik(Y)
+\cosik(Z)$, 
\item $\cosik(\{n\})=0$ for $n\in \omega$.
\end{itemize}
(2) ${\mathcal M}^{\full}$ is the set of $\cosik \in {\mathcal M}$ with 
domain ${\mathcal P}(\omega)$ and members are called ``finitely additive
measures (on $\omega$)''.

\noindent (3) We say ``$A$ has $\Xi$-measure $a$ (or $>a$, or whatever) if
$A\in\dom(\Xi)$ and $\Xi(A)$ is $a$ (or $>a$ or whatever).
\end{definition}

\begin{proposition}
\label{x1.2}
Let $a_{\alpha},b_{\alpha}$ $(\alpha<\alpha^{*})$ be reals, 
$0\leq a_{\alpha}\leq b_{\alpha}\leq 1$,  and let $A_{\alpha}\subseteq
\omega$ $(\alpha<\alpha^{*})$ be given.
The following conditions are equivalent:
\begin{enumerate}
\item[(A)] 
There exists $\cosik \in {\mathcal M}$ which satisfies $\cosik(A_{\alpha})
\in [a_{\alpha},b_{\alpha}]$ for $\alpha<\alpha^{*}$.
\item[(B)] For every $\varepsilon > 0, m<\omega$ and $n<\omega$, and
$\alpha_{0}<\alpha_{1},\ldots<\alpha_{n-1}<\alpha^{*}$ we can find a finite,  
non-empty $u\subseteq [m,\omega)$ such that for $l<n$
\[a_{\alpha_{l}}-\varepsilon\leq |A_{\alpha_{l}}\cap u|/|u|\leq
b_{\alpha_{l}}+\varepsilon.\] 
\item[(C)] For every real $\varepsilon>0$, $n<\omega$ and $\alpha_{0}<
\alpha_{1}, \ldots <\alpha_{n-1}<\alpha^{*}$ there are $c_{l}\in [
a_{\alpha_{l}}-\varepsilon,b_{\alpha_{l}}+\varepsilon]$ such that in the
vector space $\bbR^{n}$, $\langle c_{0},\ldots c_{n-1} \rangle$ is in the
convex hull of $\{\rho\in {}^{n\!}\{0,1\}:$ for infinitely many $m\in\omega$
we have $\forall_{l<n} [\rho(l)=1\Leftrightarrow m\in A_{\alpha_{l}}]\}$. 
\item[(D)] Like part (A) with $\Xi\in {\mathcal M}^{\full}$.
\item[(E)] like part (B) demanding $u\subseteq \omega$, $|u|\geq m$. 
\end{enumerate}
\end{proposition}

\begin{proof}
Straightforward. On (C) see the \ref{2.15new}.
\end{proof}

\begin{proposition}
\label{ext} 
1) Assume that $\cosik_{0}\in {\mathcal M}$ and for $\alpha<\alpha^{*}$,
$A_{\alpha} \subseteq \omega$ and $0\leq a_{\alpha} \leq b_{\alpha}
\leq 1$, $a_{\alpha},b_{\alpha}$ reals. The following are equivalent:
\begin{enumerate}
\item[(A)] There is $\cosik \in {\mathcal M}^{full}$ extending $\cosik_{0}$
such that $\alpha<\alpha^{*} \Rightarrow \cosik(A_{\alpha})\in [a_{\alpha},
b_{\alpha}]$.
\item[(B)] For every partition $\langle B_{0}, \ldots ,B_{m-1}\rangle$ of
$\omega$ with $B_{i}\in\dom(\cosik_{0})$ and $\varepsilon>0$, $n<\omega$
and $\alpha_{0}<\ldots<\alpha_{n-1}<\alpha^{*}$ we can find a finite set
$u\subseteq\omega$ such that $\cosik(B_{\alpha})-\varepsilon\leq |u\cap
B_{\alpha}|/|u|\leq\cosik(B_{\alpha})+\varepsilon$ and $a_{\alpha_{l}}-
\varepsilon \leq |u\cap A_{\alpha_{l}}|/|u|\leq b_{\alpha_{l}}+\varepsilon$.
\item[(C)] For every partition $\langle B_{0},\ldots ,B_{m-1}\rangle$ of
$\omega$ with $B_{i} \in \dom(\cosik_{0})$ and $\varepsilon>0$, $n<\omega$
and $\alpha_{0}< \ldots <\alpha_{n-1}<\alpha^{*}$ we can find $c_{l,k}\in
[0,1]_{\bbR}$ for $l<n,k<m$  such that  
\begin{enumerate}
\item[(a)] $\sum\limits_{k<m} c_{l, k}\in (a_{\alpha_l}-\varepsilon,
b_{\alpha_l}+\varepsilon)$ 
\item[(b)] for each $k<m$ and $s<\omega$ we can find $u\subseteq B_k$
with $\geq s$ members such that 
$$
l<n \Rightarrow c_{l, k}-\varepsilon<(|u\cap A_l \cap B_k|/|u|)\times
\Xi(B_k) < c_{l, k} +\varepsilon. 
$$
\end{enumerate}
\item[(D)] for every partition $\langle B_0, \ldots, B_{m-1}\rangle$
of $\omega$, with $B_i\in \dom(\Xi_0)$, $\varepsilon>0$, $n<\omega$,
and $\alpha_0, \ldots, 
\alpha_{n-1}<\alpha^*$ we can find $c_{l, k}\in [0, 1]_{\bbr}$ for
$l<n$, $k<m$ such that 
\begin{enumerate}
\item[(a)] $\sum_{k<m}c_{l,k} \in [a_{\alpha_{l}}-\varepsilon,b_{\alpha_{l}}
+\varepsilon]$,
\item[(b)] $\langle c_{l,k}:l<n\rangle$ is in the convex hull of 
$$
\begin{array}{ll}
\{\rho\in {}^{n}\{0,1\}:&\mbox{\rm for infinitely many }i\in B_k,\mbox{ we
have: }\\ 
\ &(\forall l<n)[\rho(l)=1 \Leftrightarrow i\in A_{\alpha_{l}}]\}.
\end{array}
$$ 
\end{enumerate}
\end{enumerate}
2) The following are sufficient conditions for (A), (B), (C), (D) above: 
\begin{enumerate}
\item[(E)] For every $\varepsilon > 0$, $A^{*} \in \dom(\cosik_{0})$
such that $\cosik_{0}(A^{*})>0$, $n<\omega$, $\alpha_{0}<\ldots  
<\alpha_{n-1}<\alpha^{*}$, we can find a finite, non-empty $u\subseteq 
A^{*}$ such that $a_{\alpha_{l}}-\varepsilon \leq
|A_{\alpha_{l}}\cap u|/|u|\leq b_{\alpha_{l}}+\varepsilon$ for $l<n$.
\item[(F)]  For every  $\varepsilon>0$, $n<\omega$, $\alpha_{0}<\alpha_{1},
\ldots <\alpha_{n-1}<\alpha^{*}$ and $A^{*}\in \dom(\cosik_{0})$ such that
$\cosik_{0}(A^{*})> 0$, the set $\prod_{l<n}[a_{\alpha_l}-\varepsilon,
b_{\alpha_l}+\varepsilon]\subseteq \bbR^{n}$ is not disjoint to the convex
hull of  
$$
\begin{array}{ll}
\{\rho \in {}^{n}\{0,1\}:&\mbox{\rm for infinitely many }m\in A^{*} \mbox{ we
have:}\\  
\ &\forall_{l<n}[\rho(l)=1 \Leftrightarrow m\in A_{\alpha_{l}}]\}.
\end{array}
$$
\end{enumerate}
3) If in addition $b_\alpha=1$ for $\alpha< \alpha^*$ {\em then} a sufficient
condition for (A) --- (E) above is 
\begin{enumerate}
\item[(G)] if $A^{*}\in \dom(\cosik_0)$ and $\cosik(A^{*})>0$ and
$n<\omega$ and $\alpha_0< \ldots < \alpha_{n-1}< \alpha^*$ then
$A^*\cap \bigcap\limits_{l< n} A_{\alpha_l}\neq \emptyset$. 
\end{enumerate}
\end{proposition}

\begin{proof} Straightforward. \end{proof}

\begin{definition}
(1) For $\cosik \in {\mathcal M}^{\full}$ and sequence $\bar{a}=\langle
a_{l}:l<\omega \rangle$ of reals in $[0,1]_{\bbr}$ (or just
$\sup\limits_{l<\omega}(|a_{l}|)<\infty$), let 
\[\begin{array}{l}
Av_{\cosik}(\bar{a})=\\
\sup\{\sum\limits_{k<k^{*}}\cosik(A_{k})\inf(\{a_{l}:l\in A_{k}\}):\langle
A_{k}:k<k^{*}\rangle\mbox{ is a partition of }\omega\}=\\
\inf\{\sum\limits_{k<k^{*}}\cosik(A_{k})\sup(\{a_{l}:l\in A_{k}\}):\langle
A_{k}:k<k^{*}\rangle \mbox{ is a partition of } \omega\}.
  \end{array}\]
(Easily proved that they are equal.) 

(2) For $\cosik\in \cM$, $A\subseteq \omega$  such that $\cosik(A)>0$ define
$\cosik_{A}(B)=\cosik(A \cap B)/\cosik(A)$. Clearly $\cosik_A\in \cM$ with
the same domain, $\cosik_A(A)=1$. If $B\subseteq\omega$ and $\cosik(B)>0$
then we let 
$$
Av_{\Xi\rest B}(\langle a_k: k\in B\rangle)=Av_{\Xi}(\langle a'_k: k<\omega
\rangle)/ \Xi(B)
$$ 
where
$a'_k=\left\{\begin{array}{ll} a_k & \mbox{ if }k\in B\\
0 & \mbox{ if }k\notin B \end{array}
\right.$
\end{definition}

\begin{proposition} 
\label{x} 
Assume that $\cosik\in {\mathcal M}^{\full}$ and $a^{i}_{l} \in [0,
1]_{\bbr}$ for $i<i^{*}<\omega$, $l<\omega$, $B\subseteq \omega$,
$\cosik(B)>0$ and $Av_{\cosik_{B}}(\langle a^{i}_{l}:l<\omega\rangle)=
b_{i}$ for $i<i^{*}$, $m^*<\omega$ and lastly $\varepsilon>0$. {\em Then}
for some finite $u\subseteq B\setminus m^*$ we have: if $i<i^{*}$, then
$b_{i}-\varepsilon<(\sum\{a^{i}_{l}:l\in u\})/|u|<b_{i}+\varepsilon$.
\end{proposition}

\begin{proof}
Let $B=\bigcup\limits_{j<j^{*}}B_{j}$ be a partition of $B$ with
$j^*<\omega$ such that for every $i<i^{*}$ we have 
\[\sum_{j<j^{*}}\sup\{a^{i}_{l}:l\in B_{j}\}\cosik(B_{j})-\sum_{j<j^{*}}
\inf\{a^{i}_{l}:l\in B_{j}\}\cosik(B_{j})<\varepsilon/2.\] 
Now choose $k^{*}$ large enough such that there are $k_j$ satisfying 
$k^{*}=\sum\limits_{j<j^{*}} k_{j}$ and $|k_{j}/k^{*}-\cosik(B_{j})/
\cosik(B)|<\varepsilon/(2j^{*})$ for $j<j^{*}$. Let $u_{j}\subseteq B_{j}
\setminus m^{*}$, $|u_{j}|=k_{j}$ for $j<j^{*}$. Now let
$u=\bigcup\limits_{j<j^{*}}u_{j}$. Now calculate: 
$$
\begin{array}{ll}
\sum\limits_{l\in u}a^{i}_{l}/|u| &= \sum\limits_{j<j^{*}}\sum\{a^{i}_{l}:
l\in u_{j}\}/|u|\leq \sum\limits_{j<j^{*}}\sup\{a^{i}_{l}:l\in B_{j}\}k_{j}
/k^{*}\\
\ & \leq\sum\limits_{j<j^{*}}\sup\{a^{i}_{l}:l\in B_{j}\}(\cosik(B_{j})/
\cosik(B)+\varepsilon/2j^{*})\\
\ & \leq b_i+\frac{\varepsilon}{2} + \frac{\varepsilon}{2} = b_i
+\varepsilon,
\end{array}
$$
$$
\begin{array}{ll}
\sum\limits_{l\in u}a^{i}_{l}/|u| & =\sum\limits_{j<j^{*}}\sum\{a^{i}_{l}:
l\in u_{j}\}/|u|\geq \sum\limits_{j<j^{*}}\inf\{a^{i}_{l}:l\in B_{j}\}
k_{j}/k^{*}\\
\ & \geq\sum\limits_{j<j^{*}}\inf\{a^{i}_{l}:l\in B_{j}\}(\cosik(B_{j})/
\cosik(B)-\varepsilon/(2j^{*}))>b_i-\varepsilon
\end{array}
$$
\end{proof}

\begin{claim}
\label{1.6}
Suppose $Q_1$, $Q_2$ are forcing notions, $\Xi_0\in \cM^{\full}$ in $V$, and
for $\ell=1,2$
\[\Vdash_{Q_\ell}\mbox{`` }\name{\Xi}_\ell\mbox{ is a finitely additive measure
extending }\Xi_0\mbox{ ''.}\]
{\em Then} 
\[\begin{array}{r}
\Vdash_{Q_1\times Q_2}\mbox{``there is a finitely additive measure
extending $\name{\Xi}_1$ and $\name{\Xi}_2$ }\ \\
\mbox{(hence $\Xi_0$)''.}\end{array}\]
\end{claim}

\begin{proof} Straightforward by \ref{x1.2} as:
\begin{enumerate}
\item[$(*)$] if $\Vdash_{Q_l}$ ``$\name{A}_l\subseteq\omega$'' and
$\Vdash_{Q_1\times Q_2}$ ``$\name{A}_1 \cap \name{A}_2$ is finite'' then
\[\begin{array}{r}
\Vdash_{Q_1\times Q_2}\mbox{`` for some $m$ and $A\subseteq \omega$, $A\in
V$ we have: }\ \\
\name{A}_1\setminus m\subseteq A,\ \ (\name{A}_2\setminus m)\cap A=\emptyset
\mbox{ ''.}\end{array}\] 
\end{enumerate}
\end{proof}

\begin{fact}
\label{1.7}
Assume $\Xi$ is a partial finitely additive measure, $\bar a^\alpha= \langle
a_k^\alpha: k< \omega\rangle$ sequence of reals for $\alpha< \alpha^*$
such that $\lim\sup\limits_k |a^\alpha_k|< \infty$ for each $\alpha$. {\em
Then\/} (B) $\Rightarrow$ (A) where
\begin{enumerate}
\item[(A)] there is $\Xi^*$, $\Xi\subseteq \Xi^*\in \cM^{\full}$ such
that $Av_{\Xi^*}(\bar a^\alpha)\geq b_\alpha$ for $\alpha< \alpha^*$.
\item[(B)] for every partition $\langle B_0, \ldots, B_{m^*-1}\rangle$
of $\omega$ with $B_m\in \dom(\Xi)$ and $\varepsilon>0$, $k^*>0$ and
$\alpha_0< \ldots< \alpha_{n^*-1}< \alpha^*$, there is a finite $u\subseteq
\omega\setminus k^*$ such that: 
\begin{enumerate}
\item[(i)] $\Xi(B_m)-\varepsilon<|u\cap B_m|/|u|<\Xi(B_m)+\varepsilon$, 
\item[(ii)] $\frac{1}{|u|} \sum\limits_{k\in u} a^{\alpha_l}_k>b_{\alpha_l}
-\varepsilon$ for $l< n^*$.
\end{enumerate}
\end{enumerate}
\end{fact}

\begin{remark}
If in (A) we demand $Av_{\Xi}(\bar a^\alpha)= b_\alpha$,\\
{\em then\/} in (B)(ii) add $\frac{1}{|u|} \sum\limits_{k\in u}
a^{\alpha_l}_k\leq b_{\alpha_l} +\varepsilon$.
\end{remark}

\section{The iteration}
Ignoring ${\rm MA}_{<\kappa}$ (which anyhow was a side issue) a quite
natural approach in order to get \ref{A} (i.e.\ $\cov(\nul)=\lambda$, say
$\lambda=\aleph_\omega$) is to use finite support iteration, $\bar Q=
\langle P_\alpha, Q_\alpha: \alpha< \alpha^*\rangle$, add in the first
$\lambda$ steps null sets $N_\alpha$ (the intension is that
$\bigcup\limits_{\alpha<\lambda} N_\alpha=\can$ in the final model), and
then iterate with $Q_\alpha$ being $\Random^{\V^{P'_\alpha}}$ where
$P'_\alpha\lesdot P_\alpha$ and $|P'_\alpha|<\lambda$. Say, for some
$A_\alpha\subseteq\alpha$ 
$$
\begin{array}{ll}
P'_\alpha=\{p\in P_\alpha: & \dom(p)\subseteq A_{\alpha}\mbox{ and this
holds for the conditions}\\
\ & \mbox{ involved in the }P_\gamma\mbox{-name for }\gamma\in\dom(p)
\mbox{ etc}\}
\end{array}
$$
(so each $Q_\alpha$ is a partial random; see Definition \ref{2.1}). If every
set of $<\lambda$ null sets from $\V^{P_{\alpha^*}}$ is included in some
$\V^{P'_\alpha}$, clearly $\V^{P_{\alpha^*}}\models\cov(\nul)\geq \lambda$;
but we need the other inequality too.

The problem is why does $\langle N_\alpha: \alpha<\lambda\rangle$ continue
to cover?  For $P_{\lambda+n}$ such that $\alpha\in [\lambda, \lambda+n)
\Rightarrow A_\alpha=\alpha$ this is very clear (we get iteration of Random
forcing) and if $\alpha\in [\lambda, \lambda+n)\Rightarrow A_\alpha\subseteq
\lambda$ this is clear (we get product). But necessarily we get a quite
chaotic sequence $\langle A_{\alpha_m}\cap \{\alpha_\ell: \ell<m\}:
m<m^*\rangle$ for some $\alpha_0<\ldots <\alpha_{m^*-1}$. More concretely this
is the problem of why there are no perfect sets of random reals (see
\ref{sufficient}) or even just no dominating reals. We need to ``let the 
partial randoms whisper secrets one to another", in other words to pass
information in some way. This is done by the finitely additive measures
$\name{\Xi}^t_\alpha$. We had tried with thinking of using
$\aleph_\varepsilon$-support (see \cite{Sh:538}), the idea is still clear in
the proof of \ref{3.2}. In this proof we start with ``no dominating reals"
for which we can just use ultrafilters (rather than finitely additive
measures). 

Let us start with a ground model $\V$ satisfying the following hypothesis:

\begin{hypothesis}
\label{2.0}
\begin{enumerate}
\item[(a)] $\lambda =\sum\limits_{\zeta<\delta(*)}\lambda_{\zeta}$, 
$\aleph_{0}<\kappa=\cf(\kappa)$, $\kappa<\lambda_{\zeta}<\lambda_{\gamma}$
for $\zeta <\gamma<\delta(*)$, $\chi^{\lambda}=\chi$, $2^{\kappa}= \chi$ and
$\zeta<\delta(*)\Rightarrow (\lambda_\zeta)^{\aleph_0}<\lambda$,
\item[(b)] we have one of the following\footnote{actually, any ordinal
$\alpha^*$ of cardinality $\chi$, divisible by $\chi$ and of cofinality
$>\lambda$ is O.K.}:
\begin{enumerate}
\item[($\alpha$)] $\cf(\chi)>\lambda$, the length of the final iteration is
$\chi$, 
\item[($\beta$)] length of the final iteration is $\chi\times\chi\times
\lambda^+$.
\end{enumerate}
\end{enumerate}
\end{hypothesis}
We speak mainly on $(\alpha)$. In case $(\beta)$ we should be careful to have
no repetitions in $\bar{\eta}=\langle\eta_{\alpha}:\alpha<\delta^{*}\rangle$
(see below) or $\langle \eta_{\alpha}/\approx_{\kappa}:\alpha<\delta^{*}
\rangle$ with no repetitions, where $\eta \approx_{\kappa}\nu$ iff $\eta,\nu
\in {}^{\kappa}2$ and $|\{i<\kappa:\eta(i)\neq\nu(i)\}|<\kappa$.

The reader may choose to restrict himself and start with $\V$ satisfying:
GCH, $\lambda=\aleph_\omega$, $\delta(*)=\omega$, $\lambda_{n}=\aleph_{n(*)
+n}$, $\kappa=\aleph_{n(*)} > \aleph_{1}$ and $\chi=\aleph_{\omega+1}$. Now
add $\aleph_{\omega +1}$ generic subsets of $\kappa$, i.e., force with a
product of $\chi$ copies of $(^{\kappa >}2, \vartriangleleft)$ with support
$<\kappa$. This model satisfies the hypothesis. 

We intend to define a forcing $P$ such that 
$$
\V^{P} \models 2^{\aleph_{0}}=\chi\;+\;\cov(\nul)=\lambda\;+\;{\rm
MA}_{<\kappa}. 
$$

\begin{definition}
\label{2.1}
1) ${\mathcal K}$ is the family of sequences 
$$
\bar{Q}=(P_{\alpha},\name{Q}_{\alpha},A_{\alpha},\mu_\alpha,\name{\tau}_{
\alpha}: \alpha<\alpha^*)
$$ 
satisfying: 
\begin{enumerate}
\item[(A)] $(P_{\alpha},\name{Q}_{\alpha}:\alpha<\alpha^*)$ is a finite
support iteration of c.c.c.\ forcing notions, we call $\alpha^*=\lgx(
\bar{Q})$ (the length of $\bar Q$), $P_{\alpha^*}$ is the limit,
\item[(B)] $\name{\tau}_{\alpha}\subseteq \mu_{\alpha}<\kappa$ is the
generic of $\name{Q}_{\alpha}$, (i.e.\ over $\V^{P_{\alpha}}$ from
$G_{\name{Q}_{\alpha}}$ we can compute $\name{\tau}_{\alpha}$ and vice
versa), 
\item[(C)] $A_{\alpha}\subseteq \alpha$ (for proving theorem \ref{A} we use
$|A_{\alpha}|<\lambda$),
\item[(D)] $\name{Q}_{\alpha}$ is a 
$P_\alpha$-name of a c.c.c. forcing notion 
but computable from $\langle \name{\tau}_\gamma[\name{G}_{P_\alpha}]:
\gamma\in A_\alpha\rangle$; in particular it belongs to $\V_\alpha=\V[
\langle \name{\tau}_\gamma[G_{P_\alpha}]:\gamma\in A_\alpha\rangle]$. 
\item[(E)] $\alpha^*\geq \lambda$ and for $\alpha<\lambda$ we have
$Q_{\alpha}=(^{\omega>}2,\vartriangleleft)$ (the Cohen forcing) and
$\mu_\alpha=\aleph_0$ (well, identifies $^{\omega>}2$ with $\omega$).
\item[(F)] For each $\alpha<\alpha^*$ one of the following holds, (and the
case is determined in $\V$):
\begin{enumerate}
\item[$(\alpha)$] $|\name{Q}_{\alpha}|<\kappa$, $|A_\alpha|<\kappa$ and
(just for notational simplicity) the set of elements of $\name{Q}_\alpha$ is
$\mu_\alpha<\kappa$ (but the order not necessarily the order of the
ordinals) and $\name{Q}_\alpha$ is seperative (i.e.\ $\zeta\forces
\xi\in G_{\name{Q}_{\alpha}} \Leftrightarrow \name{Q}_\alpha\models\xi\leq
\zeta$),
\item[$(\beta)$] essentially $\name{Q}_{\alpha}=\Random^{\V_{\alpha}}$ \\
($=\{r\in \V_{\alpha}:r\subseteq \fs, \mbox{ perfect tree, }\Leb(\lim(r))>
0\}$)\\
and $|A_\alpha|\geq \kappa$;
\end{enumerate}
but for simplicity $\name{Q}_\alpha=\Random^{A_\alpha, \bar Q\restriction
\alpha}$
where for $A\subseteq\lgx(\bar Q)$, 
$$
\begin{array}{ll}
\Random^{A, \bar Q}=\{p:&\mbox{there is (in $V$) a Borel function
}\cB=\cB(x_0, x_1, \ldots),\\
\ &\mbox{with variables ranging on $\{$true, false$\}$ and }\\
\ &\mbox{range perfect subtrees }r \mbox{ of }{}^{\omega>}2\mbox{ with
}\Leb(\lim r)>0,\\
\ &\mbox{such that }(\forall \eta\in r)[\Leb\lim r^{[\eta]}>0],\\
\ &\mbox{recalling }r^{[\eta]}= \{\nu \in r: \nu \XXeq \eta \vee \eta \XXeq
\nu\},\ \mbox{ and}\\ 
\ &\mbox{there are pairs }(\gamma_\ell,\zeta_\ell)\mbox{ for }\ell\!<\!\omega,
\gamma_\ell\in A\mbox{ and }\zeta_\ell<\mu_{\gamma_\ell},\\
\ &\mbox{such that }p=\cB(\ldots, \mbox{truth value}(\zeta_\ell\in
\name{\tau}_{\gamma_\ell}), \ldots)_{\ell<\omega}\}
\end{array}
$$
(in other notation, $p=\cB(\mbox{truth value}(\zeta_\ell\in\name{\tau}_{
\gamma_\ell}):\ell<\omega)$);\\
in this case we let $\supx(p)=\{\gamma_\ell: \ell<\omega\}$.\\
In this case $\mu_\alpha=\omega$ and $\name{\tau}_\alpha$ is the random real,
i.e. 
$$\name{\tau}_\alpha(n)=\ell \Leftrightarrow (\exists \eta\in {}^n 2)(^{
\omega>}2)^{[\eta\conc\langle \ell\rangle]}\in\name{G}_{Q_\alpha}.
$$
\end{enumerate}
2) Let 
$$
\begin{array}{ll}
P'_\alpha=\{p\in P_\alpha:&\mbox{for every }\gamma\in \dom(p), \mbox{if }
|A_\gamma|<\kappa\mbox{ then}\\
\ &p(\gamma)\mbox{ is an ordinal }<\mu_\gamma\mbox{ (not just a
$P_\gamma$--name) and}\\ 
\ &\mbox{if }|A_\gamma|\geq \kappa\mbox{ then }p(\gamma)\mbox{ has the form
mentioned in clause}\\
\ &\mbox{(F)$(\beta)$ above (and not just a $P_\gamma$--name of such object)
}\}
\end{array}
$$
(this is a dense subset of $P_\alpha$).

3) For $A\subseteq \alpha$ let 
$$
P'_A=\{p\in P'_\alpha: \dom(p)\subseteq A \mbox{ and }\gamma\in \dom(p)
\Rightarrow \supx(p(\gamma))\subseteq A\}
$$ 
\end{definition}

\begin{fact}\label{notl}
Suppose $\bar{Q} \in {\mathcal K}$ with $\lgx(\bar{Q})=\alpha^*$.
\begin{enumerate}
\item For $\alpha\leq\alpha^*$, $P'_\alpha$ is a dense subset of $P_\alpha$
and $P_\alpha$ satisfies the c.c.c.
\item If
\begin{enumerate}
\item[(a)]  $\cf(\alpha^*)>\lambda$,
\item[(b)] for every $A\subseteq \alpha^*$,
if $|A|<\lambda$, {\em then} there is
$\beta<\alpha^*$ such that $A\subseteq A_{\beta}$ (and $|A_\beta|\geq\kappa$).
\end{enumerate}
{\em Then}, in the extension, $\can$ is not the union of $<\lambda$ null
sets. 
\item In $\V^{P_\alpha}$, from $\name{\tau}_\alpha[G_{Q_\alpha}]$ we can
reconstruct $G_{Q_\alpha}$ and vice versa. From $\langle \name{\tau}_\gamma:
\gamma<\alpha\rangle[G_{P_\alpha}]$ we can reconstruct $G_{P_\alpha}$ and vice
versa. So $\V^{P_\alpha}=\V[\langle\name{\tau}_\beta:\beta<\alpha\rangle]$.
\item If $\mu<\lambda$, and $\name{X}$ is a $P_{\alpha^*}$-name of a subset of
$\mu$, {\em then\/} there is a set $A\subseteq\alpha^*$ such that $|A|\leq
\mu$ and $\forces_{P_{\alpha^*}}\mbox{``}\name{X}\in\V[\langle
\name{\tau}_\gamma:\gamma\in A\rangle]\mbox{''}$. Moreover for each
$\zeta<\mu$ there is in $\V$ a Borel function $\cB(x_0,\ldots,x_n,
\ldots)_{n<\omega}$ with domain and range the set $\{\true,\false\}$ and
$\gamma_\ell\in A$, $\zeta_\ell<\mu_{\gamma_\ell}$ for $\ell<\omega$ such
that 
$$
\forces_{P_{\alpha^*}}\mbox{``}\zeta\in\name{X}\mbox{ iff }\true=\cB(\ldots,
\mbox{``truth value of }\zeta_\ell\in\name{\tau}_{\gamma_\ell}[
G_{Q_{\gamma_\ell}}]\mbox{''},\ldots)_{\ell<\omega}\mbox{''}
$$
\item For $\bar Q\in {\mathcal K}$ and $A\subseteq\alpha^*$, every real in
$\V[\langle\name{\tau}_\gamma:\gamma\in A\rangle]$ has the form mentioned in
clause (F)$(\beta)$ of\ref{2.1}(1).
\item If condition (c) below holds {\em then} $\V^{P_{\alpha^*}}\models {\rm
MA}_{<\kappa}$, where
\begin{enumerate}
\item[(c)] if $\name{Q}$ is a $P_{\alpha^*}$-name of a c.c.c.\ forcing
notion with set of elements $\mu<\kappa$ {\em then} for some $\alpha<
\alpha^*$, $\name{Q}$ is a $P_\alpha$-name $\mu_\alpha=\mu$ and
$\forces_{P_\alpha}\mbox{``}\name{Q}=\name{Q}_\alpha\mbox{''}$.
\end{enumerate}
\item if $|A_\beta|\geq\kappa$ {\em then} $Q_\beta$ is actually
$\Random^{\V^{P_{A_\beta}}}$.
\end{enumerate}
\end{fact}

\begin{proof} 2) Easy using parts (3) --- (7). Note that for any $\beta<
\alpha^*$ satisfying $|A_\beta|\geq\kappa$ the null sets from $\V^{[\langle
\name{\tau}_\gamma:\gamma\in A_{\beta}\rangle]}$ do not cover $\can$ in
$\V^{P_{\alpha^*}}$ as we have random reals over $\V^{[\langle\name{
\tau}_\gamma:\gamma\in A_{\beta}\rangle]}$. So, by clause (b) of the
assumption, it is enough to note that if $\name{y}$ is a
$P_{\alpha^*}$--name of a member of $\can$, then there is a countable $A
\subseteq\alpha^*$, such that $\name{y}[\name{G}]\in\V[\langle\name{\tau}_{
\beta}:\beta\in A\rangle]$. This follows by part (4).

\noindent 3) By induction on $\alpha$.

\noindent 4) Let $\chi^*$ be such that $\{\bar Q,\lambda\}\in\cH(\chi^*)$,
and let $\zeta<\mu$; let $M$ be a countable elementary submodel of
$(\cH(\chi^*),\in, <^*_{\chi^*})$ to which $\{\bar Q,\lambda,\kappa,\mu,
\name{X},\zeta\}$ belongs, so $\forces_{P_{\alpha^*}}\mbox{``}
M[\name{G}_{P_{\alpha^*}}]\cap\cH(\chi^*)=M\mbox{''}$. Hence by \ref{notl}(3)
(i.e.\ as $\V^{P_\alpha}=\V[\langle\name{\tau}_\beta:\beta<\alpha\rangle]$)
we have $M[\name{G}_{P_{\alpha^*}}]=M[\langle\name{\tau}_i: i\in\alpha^*\cap
M\rangle]$ and the conclusion should be clear.

\noindent 5) By \ref{notl}(4).

\noindent 6) Straight.

\noindent 7) Check.
\end{proof}

\begin{definition}
\label{2.3}
(1) Suppose that $\bar{a}=\langle a_{l}:l<\omega\rangle$ and $\langle
n_l: l<\omega\rangle$ are such that:
\begin{enumerate}
\item[(a)] $a_{l}\subseteq {}^{n_{l}}2$,
\item[(b)] $n_{l}<n_{l+1}<\omega$ for $l<\omega$,
\item[(c)] $|a_{l}|/2^{n_{l}}>1-1/ 10^{l}$.
\end{enumerate}
Let $N[\bar{a}]=:\{\eta\in\can : (\exists^{\infty} l)(\forall\nu\in a_{l})
\nu\not\vartriangleleft\eta\}$.

\noindent (2) For $\bar{a}$ as above and $n\in\omega$, let 
$$
\tree_{n}(\bar{a})=
\{\nu\in\fs : n_{l}>n\Rightarrow\nu\restriction n_{l}\in a_{l}\}.
$$
\end{definition}

\medskip
It is well known that for $\bar{a}$ as above the set $N[\bar{a}]$ is null
(and $N[\bar a]= {}^\omega 2\setminus\bigcup\limits_{n<\omega}\lim
\tree_n(\bar a)$).

\begin{definition}
\label{2.4}
For $\alpha<\lambda$ we identify $Q_{\alpha}$ (the Cohen forcing) with:
$$
\{\langle(n_{l},a_{l}):l<k\rangle: k<\omega, n_{l}<n_{l+1}<\omega,\;
a_{l}\subseteq {}^{n_{l}}2,\; |a_{l}|/2^{n_{l}}>1-1/10^{l}\},
$$ 
ordered by the end extension. If $G_{Q_\alpha}$ is $Q_{\alpha}$-generic, let
$\langle(\name{\bar{a}}^{\alpha}_\ell,\name{n}^\alpha_\ell):\ell<\omega\rangle
[G_{Q_\alpha}]$ be the $\omega$--sequence such that every $p\in
G_{Q_\alpha}$ is an initial segment of it. So we have defined the
$Q_\alpha$--name $\bar a^\alpha=\langle\name{a}^\alpha_\ell:\ell<\omega
\rangle$ and similarly $\langle\name{n}_\ell^\alpha: \ell<\omega\rangle$.
Let $\name{N}_{\alpha}=N[\name{\bar{a}}^{\alpha}]$.
\end{definition}

Our aim is to prove that 

\begin{definition}
\label{2.5A}
For $\bar Q\in\cK$ with $\alpha^*=\lgx(\bar Q)$ let:
\begin{enumerate}
\item[$(*)_{\bar{Q}}$] $\langle N_{\alpha}:\alpha <\lambda\rangle$ cover
$\can$ in $\V^{P_{\alpha^*}}$, where $P_{\alpha^*}=\Lim(\bar{Q})$.
\end{enumerate}
\end{definition}

We eventually shall prove it not for every $\bar Q$, but for enough $\bar
Q$'s (basically asking the $A_\alpha$ of cardinality $\geq\kappa$ to be
closed enough).

\begin{lemma}
\label{sufficient}
For $\bar Q\in\cK$ with $\gamma=\lgx(\bar Q)$, a sufficient condition for
$(*)_{\bar{Q}}$ is: 
\begin{enumerate}
\item[$(**)_{\bar{Q}}$] In $\V^{P_{\gamma}}$: there is no perfect tree
$T\subseteq\fs$ and $E\in [\lambda]^{\kappa^{+}}$ such that, for some
$n<\omega$, $T\subseteq\tree_{n}[\bar{a}^{\alpha}]$ for all $\alpha\in E$.
\end{enumerate}
\end{lemma}

\begin{proof}
By induction on $\gamma\geq\lambda$. For $\gamma=\lambda$, trivial by
properties of the Cohen forcing.

Suppose $\gamma>\lambda$ limit. Assume toward contradiction that
$$
p\forces_{P_{\gamma}}\mbox{`` }\name{\eta}\not\in
\bigcup_{\alpha <\lambda}N[\name{\bar{a}}^{\alpha}]\mbox{ ''}.
$$
W.l.o.g. 
$$
p\forces_{P_{\gamma}}\mbox{`` }\name{\eta}\not\in {\V^{P_{\beta}}}\mbox{ ''}
$$ 
for every $\beta<\gamma$, hence by properties of FS iteration of c.c.c.\
forcing notions $\cf(\gamma)=\aleph_0$. So for each $\alpha <\lambda$ there
are $p_{\alpha}$, $m_{\alpha}$ such that 
$$
p\leq p_{\alpha}\in P_{\gamma},\qquad p_{\alpha}\forces\mbox{`` }\name{\eta}
\in\lim(\tree_{m_{\alpha}}(\name{\bar{a}}^{\alpha}))\mbox{ ''}.
$$
Note that (by properties of c.c.c.\ forcing notions) $\langle\{\alpha<
\lambda: p_{\alpha}\in P_{\beta}\}:\beta<\gamma\rangle$ is an increasing
sequence of subsets of $\lambda$ of length $\gamma$, so for some $\gamma_{1}
<\gamma$ there is $E\in [\lambda]^{\kappa^{+}}$ such that $p_{\alpha}\in
P_{\gamma_{1}}$ for every $\alpha\in E$ and w.l.o.g. $m_{\alpha}=m$ for
$\alpha\in E$. Note that for all but $<\kappa^{+}$ of the ordinals
$\alpha\in E$ we have 
$$
p_{\alpha}\forces{\rm ``}|\{\beta\in E: p_{\beta}\in G_{P_{\gamma_{1}}}\}|
=\kappa^{+}\mbox{\rm ''}.
$$ 
Fix such $\alpha$, let $G_{P_{\gamma_{1}}}$ be $P_{\gamma_{1}}$-generic over
$\V$ subset of $P_{\gamma_1}$ to which $p_\alpha$ belongs. Now in
$\V[G_{P_{\gamma_1}}]$ let $E'=\{\beta\in E:p_{\beta}\in G_{P_{\gamma_{1}
}}\}$ so $|E'|=\kappa^+$. Let $T^{*}=\bigcap_{\beta\in E'}\tree_{m}(
\bar{a}^{\beta})$. Note that, in $\V^{P_{\gamma_{1}}}$, $T^{*}$ is a subtree
of $\fs$ and by $(**)$, $T^{*}$ contains no perfect subtree. Hence
$\lim(T^{*})$ is countable, so absolute. But $p_{\alpha}\forces_{P_\gamma}
\mbox{`` }\name{\eta}\in\lim(T^{*})\mbox{ ''}$, so $p_{\alpha}\forces\mbox{
``}\name{\eta}\in {V^{P_{\gamma_{1}}}}\mbox{ ''}$, a contradiction.

Assume now that $\gamma=\beta+1>\lambda$ and work in $V^{P_{\beta}}$. Choose
$p,p_{\alpha}\in Q_{\beta}$ as before. Note that $Q_{\beta}$ has a dense
subset of cardinality $<\lambda$, so there is some $q\in Q_{\beta}$ and $m$
such that $E=\{\alpha <\lambda:m_{\alpha}=n,\ p_{\alpha}\leq q\}$ has
cardinality $\geq\kappa^{+}$. Continue as above.

As we have covered the cases $\gamma=\lambda$, $\gamma>\lambda$ limit and
$\gamma>\lambda$ successor, we have finished the proof.
\end{proof}

\begin{discussion}
\label{2.6A}
Note that by Lemma \ref{sufficient} and Fact \ref{notl} it is enough to show
that there is $\bar{Q}\in {\mathcal K}$ such that $\alpha^*=\lgx(\bar{Q})$
(where $\alpha^*$ is chosen as the length of the final iteration from
\ref{2.0} clause (b)), satisfying clauses (a)+(b)+(c) of Fact \ref{notl}(2)
+ \ref{notl}(6) and $(**)_{\bar Q}$. To prove the latter we need to impose
more restrictions on the iteration.  
\end{discussion}

\begin{definition} \label{T} 
$\cT$, the set of blueprints, is the set of tuples
$$
t=(w^{t},\n^{t},\m^t,\bar{\eta}^{t},h_{0}^{t},h_{1}^{t},h^{t}_{2},\bar{n}^t)
$$
where: 
\begin{enumerate}
\item[(a)] $w^{t}\in [\kappa]^{\aleph_{0}}$,
\item[(b)] $0<\n^{t}<\omega$, $\m^t\leq\n^t$
\item[(c)] $\bar{\eta}^{t}=
\langle\eta^{t}_{\n,k}:\n<\n^{t},k<\omega\rangle$, $\eta^{t}_{\n, k}\in 
{}^{w^{t}\!}2$ for $\n<\n^{t}, k<\omega$,
\item[(d)] $h^{t}_{0}$ is a partial function from $[0,\n^{t} )$ to $\kom$,
its domain includes the set $\{0,\ldots,\m^t-1\}$ (here we consider members
of $Q_{\alpha}$ (for $\alpha <\lambda$) as integers\footnote{ actually the
case where each $h^t_0(\n)$ is a constant function from $\omega$ to $\kappa$
suffices, and so $\kappa<\lambda$ suffices instead $\kappa^{\aleph_0}<
\lambda$}),
\item[(e)] $h^{t}_{1}$ is a partial function from $[0,\n^{t})$ to $(0,1)_{
\bbQ}$ (rationals), but for $\n\in [0,\n^t)\setminus\dom(h_1^t)$ we
stipulate $h^t_1(\n)=0$ and we assume $\sum\limits_{\n<\n^t}
\sqrt{h^t_1(\n)}< 1/10$.  
\item[(f)] $h^{t}_{2}$ is a partial function from $[0,\n^{t})$ to $\fs$,
\item[(g)] $\dom(h^{t}_{0})$, $\dom(h^{t}_{1})$ are disjoint with union
$[0,\n^{t})$,
\item[(h)] $\dom(h^{t}_{2})=\dom(h^{t}_{1})$,
\item[(i)] $\eta^{t}_{\n_{1},k_{1}}=\eta^{t}_{\n_{2},k_{2}}\Rightarrow
\n_{1}=\n_{2}$,
\item[(j)] for each $\n<\n^t$ we have: $\langle\eta^{t}_{\n, k}: k<\omega
\rangle$ is constant or with no repetitions; if it is constant and $\n\in
\dom(h^t_0)$ then $h^t_0$ is constant.
\item[(k)] $\bar n^t=\langle n^t_k: k<\omega\rangle$ where $n^t_0=0$, $n^t_k<
n^t_{k+1}<\omega$ and the sequence $\langle n^t_{k+1}-n^t_k: k<\omega\rangle$
goes to infinity. Let for $\ell<\omega$ and such $\bar n$, $k_{\bar n}(\ell)=
k(\ell,\bar n)$ be the unique $k$ such that $n_k\leq\ell< n_{k+1}$.
\end{enumerate}
\end{definition}

\begin{discussion}
\label{2.9new}
The definitions of a blueprint $t\in\cT$ (in Definition \ref{T}) and of
iterations $\bar Q\in\cK^3$ (defined in Definition \ref{K3} clause (c)
below; the reader may first read it) contain the main idea of the proof, so
though they have many clauses, the reader is advised to try to understand
them. 

In order to prove $(**)_{\bar{Q}}$ we will show in $\V^{P_{\alpha^*}}$ that
if $E\in[\lambda]^{\kappa^+}$, and $n<\omega$, then $\bigcap\limits_{\alpha
\in E}\tree_n(\name{\bar a}^\alpha)$ is a tree with finitely many branches. 
So let $p$ be given, let $p\leq p_\zeta\forces\mbox{`` }\beta_\zeta\in
\name{E}$ '' for $\zeta<\kappa^+$, $\beta_\zeta\notin\{\beta_\xi:\xi<\zeta
\}$, we can assume $p_\zeta$ is in some pregiven dense set, and $\langle
p_\zeta:\zeta<\kappa^+\rangle$ form a $\Delta$-system (with some more
``thinning" demands), $\dom(p_\zeta)=\{\alpha_{\n,\zeta}:\n<\n^*\}$,
$\alpha_{\n,\zeta}$ is increasing with $n$, and $\alpha_{\n,\zeta}<\lambda$
iff $\n<\m^*$. Let $p^\prime_\zeta$ be $p_\zeta$ when $p_\zeta(\alpha_{\n,
\zeta})$ is increased a little, as described below. 

It suffices to find $p^*\geq p$ such that $p^*\forces\mbox{`` }\name{A}=:
\{\zeta<\omega: p^\prime_\zeta\in\name{G}\}$ is large enough such that
$\bigcap\limits_{\zeta\in\name{A}}\tree_m(\bar a^{\beta_\zeta})$ has only
finitely many branches''. 

Because of ``communication problems'' the ``large enough'' is interpreted as
of $\name{\Xi}^t_\alpha$-measure (again defined in \ref{K3} below).

The natural numbers $\n<\n^*$ such that $Q_{\alpha_{\n,\zeta}}$ is a forcing
notion of cardinality $<\kappa$, do not cause problems, as $h^t_0(\n)$
tells us exactly what the condition $p_\zeta(\alpha_{\n,\zeta})$ is. Still
there are many cases of such $\langle p_\zeta:\zeta<\omega\rangle$ which
fall into the same $t$; we possibly will get contradictory demands if
$\alpha_{\n_1,\zeta_1}=\alpha_{\n_2,\zeta_2}$, $\n_1\neq\n_2$. But the
$w^t$, $\bar \eta^t$ are exactly built to make this case not to happen. That
is, we have to assume $2^\kappa=\chi$ ($=|\alpha^*|$) in order to be able
for our iteration $\langle P_\alpha,\name{Q}_\alpha:\alpha<\alpha^*\rangle$
to choose $\langle\eta_\alpha:\alpha<\alpha^*\rangle$, $\eta_\alpha\in
{}^\kappa2$ with no repetitions, so that if $v\subseteq\chi$, $|v|\leq
\aleph_0$ (e.g.\ $v=\{\alpha_{\n, \zeta}: \n<\n^t, \zeta<\omega\}$) then for
some $w=w^t\in [\kappa]^{\aleph_0}$ we have $\langle \eta_\alpha\restriction
w: \alpha\in v\rangle$ is with no repetitions. 

So the blueprint $t$ describes such situation, giving as much information as
we can, as long as the number of blueprints is not too large, $\kappa^{
\aleph_0}=\kappa$ in our case.

If $Q_{\alpha_{\n, \zeta}}$ is a partial random, we may get many candidates
for $p_\zeta(\alpha_{\n, \zeta})\in \Random$ and they are not all the same
ones. We want that in many cases they will be in the generic set. Well, we
can (using $h^t_1(\n)$, $h_2^t(\n)$) know that in some interval
$(^\omega2)^{[h^t_2(\n)]}$ the set $\lim p_\zeta(\alpha_{\n, \zeta})$ is
large, say of relative measure $\geq 1-h^t_1(\n)$, and we could have chosen
the $p_\zeta$'s such that $\langle h^t_1(\n): \n<\n^t\rangle$ is small
enough, still the number of candidates is not bounded by $1/h^t_1(\n)$. Here
taking limit by ultrafilters is not good enough, but using finitely additive
measures is. 

Well, we have explained $w^t$, $\bar \eta^t$, $h_0^t$, $h^t_1$, $h^t_2$, but
what about the $\bar n^t=\langle n^t_k: k<\omega\rangle$? In the end (in \S
3) the specific demand on $\{\zeta: p'_\zeta\in G\}$ being large, is that
for infinitely many $k<\omega$,
$$|\{\ell: n^t_k\leq \ell< n^t_{k+1}\mbox{ and } p_\ell\in G\}| /
(n^t_{k+1} - n^t_k)
$$
is large, the $n^t_k$ will be chosen such that it is increasing fast enough
and $\langle p'_\ell(\beta_\ell): \ell\in [n^t_k,n^t_{k+1})\rangle$ will be
chosen such that for each $\varepsilon>0$ for some $s<\omega$, for $k$ large
enough if the above fraction is $\geq \varepsilon$ then essentially $^k2
\cap\{\tree_m(\bar a^{\beta_\ell}): n^t_k\leq \ell < n^t_{k+1}\mbox{ and }
p'_\ell\in G\}$ has $\leq s$ members, this suffices.

What is our plan? We define $\cK^3$, the class of suitable expanded
iterations $\bar Q$ by choice of $\eta_\alpha$ (for $\alpha< \lgx (\bar Q)$)
and names for finitely additive measures $\name{\Xi}^t_\alpha$ satisfying the
demands natural in this context. You may wonder why we use $\Xi$-averages;
this is like integral or expected value, and so ``behave nicely'' making the
``probability computations'' simpler.  Then we show that we can find $\bar
Q\in \cK^3$ in which all obligations toward ``$\cov(\nul)\geq \lambda$'' and
${\rm MA}_{<\kappa}$ hold.

The main point of \S 3 will be that we can carry the argument of ``for some
$p^*$ we have $p^*\Vdash\{\ell<\omega: p'_\ell\in G\}$ is large'' and why it
gives $n<\omega\ \&\ E\in [\chi]^{\kappa^+} \Rightarrow\bigcap\limits_{\zeta
\in A}\tree_n (\bar a^\zeta)$ has finitely many branches, thus proving
theorem \ref{A}. 

The reader may wonder how much the $\name{\Xi}^t_\alpha$ are actually
needed. As explained above they are just a transparent way to express
the property; this will be utilized in \cite{Sh:619}.
\end{discussion}

\begin{definition} \label{K3}
${\mathcal K}^{3}$ is the class of sequences
$$
\bar{Q}=\langle
P_{\alpha},\name{Q}_{\beta},A_{\beta}, \mu_\beta, \name{\tau}_{\beta},
\eta_{\beta},(\name{\Xi}_{\alpha}^{t})_{t\in \cT}:\alpha \leq \alpha^{*},
\beta <\alpha^{*}\rangle
$$ 
(we write $\alpha^{*}=\lgx(\bar{Q})$), such that:
\begin{enumerate}
\item[(a)] $\langle P_{\alpha},\name{Q}_{\beta}, A_{\beta}, \mu_\beta, 
\name{\tau}_{\beta}: 
\alpha\leq \alpha^{*},\beta<\alpha^{*}\rangle$ is in ${\mathcal K}$,
\item[(b)] $\eta_{\beta}\in \tka$ and $\forall {\beta < \alpha <\alpha^{*}}
[\eta_{\alpha}\not =\eta_{\beta}]$,
\item[(c)] $\cT$ is from Definition \ref{T}, and $\name{\Xi}_{\alpha}^{t} $
is a $P_{\alpha}$-name of a finitely additive measure on
$\omega$ (in $\V^{P_\alpha}$), increasing with $\alpha$, 
\item[(d)] We say that $\bar{\alpha}=
\langle \alpha_{l}:l<\omega \rangle$ satisfies $(t, \n)$ (for $\bar{Q}$) if:
\begin{itemize}
\item $\langle \alpha_l: l<\omega\rangle\in \V$ (of course),
\item $t\in \cT$, $\n<\n^t$,
\item $\alpha_{l}\leq \alpha_{l+1} < \alpha^{*}$,
\item $\n<\m^{t} \Leftrightarrow (\forall l)(\alpha_{l}<\lambda)
\Leftrightarrow (\exists k)(\alpha_{k}<\lambda)$,
\item $\eta^{t}_{\n, l} \subseteq \eta_{\alpha_{l}}$ (as functions),
\item if $\n \in \dom(h^{t}_{0})$, {\em then} $\mu_{\alpha_l}<\kappa$ and 
$$
\Vdash_{P_{\alpha_{l}}}\mbox{``}|Q_{\alpha_{l}}|<\kappa\mbox{ and }
(h^{t}_{0}(\n))(l)\in\name{Q}_{\alpha_{l}}\ \mbox{ i.e. }(h^t_0(\n))(l)<
\mu_{\alpha_l}\mbox{''},
$$
\item if $\n \in \dom(h^{t}_{1})$, {\em then} $\mu_{\alpha_l}\geq \kappa$ so
$ \Vdash_{P_{\alpha_{l}}}`` Q_{\alpha_{l}} \mbox{ \rm has cardinality } \geq
\kappa\mbox{\rm ''}$,  (so it is a partial random),
\item if $\langle \eta^{t}_{\n, k}: k\in \omega \rangle$ is constant, then 
$\forall {l}[\alpha_{l}=\alpha_{0}]$, if $\langle \eta^{t}_{\n, k}: k <
\omega \rangle$ is not constant, then $\forall l [\alpha_{l}<\alpha_{l+1}]$,
\end{itemize}
\item[(e)] if $\bar{\alpha}=\langle \alpha_{l}:l<\omega \rangle$ satisfies
$(t,\n)$ for $\bar{Q}$, $\bigwedge\limits_{l<\omega}(\alpha_l<\alpha_{l+1})$, 
$\n\in \dom(h^{t}_{0})$ then
$$
\begin{array}{ll}
\Vdash_{P_{\alpha^{*}}}&\mbox{ ``the following set has }\name{\Xi}^t_{
\alpha^*}\mbox{- measure }1:\\
\ & \{k<\omega: \mbox{ if }l\in [n^t_k,n^t_{k+1})\mbox{ then }(h^t_0(\n))(l)
\in\name{G}_{Q_{\alpha_l}}\}\mbox{''},
\end{array}
$$
\item[(f)] {\em if} $\bar{\alpha}=\langle \alpha_{l}:l<\omega \rangle$
satisfies $(t,\n)$ for $\bar{Q}$, $\n \in \dom(h^{t}_{1})$,
$\forall {l<\omega}[\alpha_{l} <\alpha_{l+1}]$, and $\bar r=\langle
\name{r}_l:l<\omega\rangle$ where for $l<\omega$, $\name{r}_{l}$ is a
$P_{\alpha_{l}}$-name of a member of $\name{Q}_{\alpha_{l}}$ such that (it
is forced that) 
\begin{enumerate}
\item[$(*)$] $1-h^{t}_{1}(\n) \leq \Leb(\{\eta\in \can:h^{t}_{2}(\n)
\vartriangleleft\eta \in \lim(\name{r}_{l})\})/2^{\lgx(h^{t}_{2}(\n))}$,
\end{enumerate}
\noindent {\em then} for each $\varepsilon>0$ we have
$$
\begin{array}{ll}
\Vdash_{P_{\alpha^{*}}}&\mbox{`` the following set has }\name{\Xi}^t_{
\alpha^*}\mbox{-measure }1:\\
\ &\ \ \{k<\omega: \mbox{in the set }\{\ell\in [n^t_k,n^t_{k+1}):
\name{r}_\ell\in G_{Q_{\alpha_\ell}}\}\mbox{ there are}\\ 
\ &\qquad\mbox{ at least }(n^t_{k+1}-n^t_k)\times (1-h^t_1(\n))\times
(1-\varepsilon)\mbox{ elements}\}\mbox{ ''} 
\end{array}
$$
\item[(g)] {\em if} $\bar{\alpha}=\langle\alpha_{l}:l<\omega\rangle$
satisfies $(t, \n)$ for $\bar{Q}$, $\n\in \dom(h^{t}_{1})$, $\forall l
[\alpha_{l}=\alpha]$, and $\name{r},\name{r}_{l}$ are $P'_{A_\alpha}$-names
of members of $Q_{\alpha}$ satisfying $(**)^{\bar{Q}}_{\name{r},\langle
\name{r}_l:l<\omega\rangle}$ (see below for the definition of $(**)$) {\em
then}  
$$
\begin{array}{ll}
\Vdash_{P_{\alpha^*}}& \mbox{`` if }\name{r}\in G_{Q_\alpha}\mbox{ then }
\ (1-h^t_1(\n)) \\
\ & \quad\leq Av_{\name{\Xi}^{t}_{\alpha^*}}(\langle |\{l\in
[n^t_k, n^t_{k+1}): \name{r}_l\in \name{G}_{Q_\alpha}\}|/(n^t_{k+1}- n^t_k):
k<\omega\rangle)\mbox{ ''},
\end{array}
$$
where
\begin{enumerate}
\item[$(**)^{\bar{Q}}_{\name{r},\langle\name{r}_l:l<\omega\rangle}$]
$\name{r}$, $\name{r_l}$ are $P'_{A_\alpha}$--names of members of
$Q_\alpha$ and, in $\V^{P_{\alpha}}$, for every $r'\in Q_{\alpha}$
satisfying $r\leq r'$ we have 
$$
Av_{\Xi^t_{\alpha}}(\langle a_k(r'): k<\omega\rangle)\geq 
(1-h^t_1(\n))
$$ 
where 
$$
\begin{array}{ll}
(\boxtimes)&\qquad a_k(r')= a_k(r',\bar r)=a_k(r',\bar{r},\bar{n}^t)\\
 &\qquad=\left(\sum\limits_{\ell\in [n^t_k, n^t_{k+1})} \frac{\Leb(
\lim(r')\cap\lim(r_\ell))}{\Leb(\lim(r'))}\right)\times \frac{1}{n^t_{k+1}-
 n^t_k} 
\end{array}
$$
\end{enumerate}
(so $a_k(r',\bar{r},\bar{n})\in [0,1]_{\mathbb R}$ is well defined for
$k<\omega$, $\bar{r}=\langle r_l:l<\omega\rangle$, $\{r,r_l\}\subseteq
\Random$, $\bar{n}=\langle n_l:l<\omega\rangle$, $n_l<n_{l+1}<\omega$),
\item[(h)] $P'_{A_{\alpha}} \lesdot P_{\alpha}$, and $\beta\in A_\alpha\ \&\
|A_\beta|<\kappa\ \Rightarrow\ A_\beta\subseteq A_\alpha$,
\item[(i)] if $\Vdash_{P_\alpha}$ ``$|\name{Q}_{\alpha}| \geq \kappa$", {\em
then} $\name{\Xi}^{t}_{\alpha} \restriction {\mathcal
P}(\omega)^{\V^{P_{A_{\alpha}}}}$ is a $P_{A_{\alpha}}$-name\footnote{Here
the secret was whispered.}, for every $t\in\cT$. 
\end{enumerate}
\end{definition}

\begin{definition}
\label{2.9}
\begin{enumerate}
\item For $\bar{Q}\in {\mathcal K}^{3}$ and
for $\alpha^* \leq \lgx(\bar{Q})$ let
\[\bar{Q}\restriction \alpha^* =\langle P_{\alpha},\name{Q}_{\beta},
A_{\beta}, \mu_\beta, \name{\tau}_{\beta},\eta_{\beta},(\Xi_{\alpha}^{t}
)_{t\in \cT}:\alpha \leq \alpha^*,\beta <\alpha^*\rangle,\] 
\item For $\bar Q^1$, $\bar Q^2\in \cK^3$ we say:
\[\bar{Q}^{1} \leq \bar{Q}^{2}\quad\mbox{ if }\ \bar{Q}^{1}=
\bar{Q}^{2}\restriction \lgx(\bar{Q}^{1}).\]
\end{enumerate}
\end{definition}

\begin{fact}
\label{2.9A}
(1) If $\bar{Q}\in {\mathcal K}^{3}$, $\alpha \leq \lgx(\bar{Q})$, 
{\em then} $\bar{Q}\restriction \alpha\in {\mathcal K}^{3}$.

\noindent (2) $({\mathcal K}^{3},\leq)$ is a partial order.

\noindent (3) If a sequence $\langle\bar{Q}^{\beta}:\beta<\delta\rangle
\subseteq {\mathcal K}^{3}$ is increasing, $\cf(\delta) > \aleph_{0}$, {\em
then\/} there is a unique $\bar{Q} \in {\mathcal K}^{3}$ which is the least
upper bound, $\lgx(\bar{Q})=\bigcup\limits_{\beta<\delta}\lgx(\bar{Q}^{
\beta})$ and $\bar{Q}^{\beta}\leq\bar{Q}$ for all $\beta<\delta$. 
\end{fact}

\begin{proof} Easy (recall that it is well known that $(\can)^{\V^{P_\delta}}
=\bigcup\limits_{\beta<\delta}(\can)^{\V^{P_\beta}}$, so
$\name{\Xi}^t_\delta=\bigcup\limits_{\beta<\delta}\Xi^t_\beta$ is a legal
choice). 
\end{proof} 

\begin{lemma}
\label{2.10}
Suppose that $\bar{Q}^{n}\in {\mathcal K}^{3}$, $\bar{Q}^{n}<\bar{Q}^{n+1}$,
$\alpha_{n}=\lgx(\bar{Q}^{n})$,
$\delta=\sup\limits_{n<\omega}(\alpha_{n})$. {\em Then} there is 
$\bar{Q}\in {\mathcal K}^{3}$ such that $\lgx(\bar{Q})=\delta$ and 
$\bar{Q}^{n}\leq \bar{Q}$ for $n\in\omega$.
\end{lemma}

\begin{proof} Note that the only problem is to define
$\name{\Xi}^{t}_{\delta}$ for $t\in \cT$, i.e., we have to extend
$\bigcup\limits_{\alpha<\delta}\name{\Xi}^{t}_{\alpha}$ so that the
following two conditions are satisfied, (they correspond to clauses (f) and
(e) of Definition \ref{K3}). 
\begin{enumerate}
\item[(a)] we are given\footnote[4]{in $\V$, so $\langle (\alpha_l,p_l):l<
\omega\rangle\in\V$, of course} $\n<\n^{t}$, $\n\in\dom(h^t_1)$,
$\langle\alpha_{l}:l<\omega\rangle$, from $\V$ of course, satisfies $(t,\n)$
for $\bar Q$ and is strictly increasing with limit $\delta$ and we are given
$\langle\name{p}_{l} :l<\omega \rangle$ such that $\Vdash_{P_{\alpha_l}}{\rm
`` }\name{p}_l\in \name{Q}_{\alpha_l}\mbox{ and }1-h^{t}_{1}(\n) \leq
\Leb(\{\eta\in \can:h^{t}_{2}(\n)\vartriangleleft\eta \in
\lim(\name{p}_{l})\})/2^{\lgx(h^{t}_{2}(\n))}$ ''. The demand is: for each
$\varepsilon>0$ we have $\Vdash_{P_\delta}\mbox{\rm `` }\name{\Xi}^t_\delta(
\name{C})=1$ '', where 
$$
\begin{array}{ll}
\name{C}=\{k<\omega:&\mbox{in the set }\{\ell: \ell\in
[n^t_{k}, n^t_{k+1})\mbox{ and }p_\ell\in \name{G}_{Q_{\alpha_\ell}}\}
\mbox{ there are}\\
\ &\mbox{at least }(n^t_{k+1}-n^t_k)\times (1-h^t_1(\n))\times
(1-\varepsilon)\mbox{ elements}\}.
\end{array}
$$
\item[(b)] If$^4$ $\n<\n^{t}$, $\n\in \dom(h^{t}_{0})$, $\langle
\alpha_\ell: \ell<\omega\rangle$ satisfies $(t, \n)$ for $\bar Q$ and is
strictly increasing with limit $\delta$, and $p_{l}\in Q_{\alpha_{l}}$,
satisfy $p_{l}=h^t_0(\n)(l)$ for $\ell<\omega$ (an ordinal
$<\mu_{\alpha_l}$), {\em then} 
$\Vdash _{P_{\delta}} ``\name{\Xi}^t_\delta (\name{C})= 1$'' where
$$
\name{C}=\{k<\omega: \mbox{ for every }l\in [n_k, n_{k+1})\mbox{ we
have }p_l\in \name{G}_{Q_{\alpha_l}}\}
$$
\end{enumerate}

As $\bigcup\limits_{\alpha<\delta}\name{\Xi}^t_\alpha$ is a
($P_\delta$-name of a) member of $\cM$, in $\V^{P_\delta}$
by \ref{ext}(3) it suffices to prove
\begin{enumerate}
\item[$(*)$] $\Vdash_{P_\delta}\mbox{``if }\name{B}\in
\bigcup\limits_{\alpha<\delta}\dom(\name{\Xi}^t_\alpha)=
\bigcup\limits_{\alpha<\delta} 
\cP(\omega)^{\V[P_\alpha]}$ and $\name{\Xi}^t_\alpha(\name{B})>0$
and $j^*<\omega$, and
$\name{C}_j$ (for $j<j^*$) are from (a), (b) above {\em then}
$\name{B}\cap \bigcap\limits_{j<j^*} \name{C}_j\neq \emptyset$''.
\end{enumerate}
Toward contradiction assume $q\in P_\delta$ force the negation so possibly
increasing $q$ we have: for some $\name{B}$ and for some $j^*<\omega$, for
each $j<j^*$ we have the $\varepsilon_j>0$, and $\n(j)<\n^t$, $\langle
\alpha^j_l: l<\omega\rangle$, $\langle p^j_l: l<\omega\rangle$ involved in
the definition of $\name{C}_j$ (in (a) or (b) above), $q$ force:
$\name{B}\in \bigcup\limits_{\alpha<\delta}\dom(\name{\Xi}^t_\alpha)=
\bigcup\limits_{\alpha<\delta}\cP(\omega)^{\V[P_\alpha]}$ and
$(\bigcup\limits_{\alpha<\delta} \name{\Xi}^t_\alpha)(\name{B}) > 0$ and
$\name{C}_j$ (for $j<j^*$) comes from (a) or (b) above, but $\name{B}\cap
\bigcap\limits_{j<j^*} \name{C}_j=\emptyset$; as we can decrease
$\varepsilon$, wlog $\varepsilon_j=\varepsilon$. Again w.l.o.g.  for some
$\alpha(*)< \delta$ we have $\name{B}\in \dom(\name{\Xi}^t_{\alpha(*)})$ is
a $P_{\alpha(*)}$-name, and $\name{C}_j$ have the $\n(j)< \n^t$, $\langle
\alpha^j_l: l<\omega\rangle$, $\langle p^j_l: l<\omega\rangle$ witnessing it
is as required in (a) or (b) above. W.l.o.g.\ $q\in P_{\alpha(*)}$. Possibly
increasing $q$ (inside $P_{\alpha(*)}$ though) we can find $k<\omega$ such
that $q\vDash$``$k\in\name{B}$'' and $\bigwedge\limits_{j<j^*}
\bigwedge\limits_{l\in [n^t_k, n^t_{k+1})}\alpha^j_l> \alpha(*)$ and
moreover such that $n^t_{k+1}-n^t_k$ is large enough compared to $1/
\varepsilon$, $j^*$ (just let $q\in G_{P_{\alpha(*)}}\subseteq
P_{\alpha(*)}$, $G_{P_{\alpha(*)}}$ generic over $\V$ and think in
$\V[G_{P_{\alpha(*)}}]$). Let $\{\alpha^j_l: j<j^*\mbox{ and }l\in [n^t_k,
n^t_{k+1})\}$ be listed as $\{\beta_m: m<m^*\}$, in increasing order (so
$\beta_0>\alpha(*)$) (possibly $\alpha^{j(1)}_{l(1)}=\alpha^{j(2)}_{l(2)}\
\&\ (j(1), l(1))\neq (j(2), l(2))$). Now we choose by induction on $m\leq
m^*$ a condition $q_m\in P_{\beta_m}$ above $q$, increasing with $m$, where
we stipulate $\beta_{m^*}=\delta$.

During this definition we ``throw a dice'' and prove that the probability of
success (i.e.\ $q_{m^*}\Vdash ``k\in \name{C}_j$'' for $j<j^*$) is positive,
so there is $q_{m^*}$ as desired hence we get the desired contradiction.
\medskip

\noindent {\em Case A:} $m=0$

Let $q_0=q$
\medskip

\noindent {\em Case B:} $m+1$, and for some $\n< \n^t$, we have $\n\in
\dom(h^t_0)$ and $\zeta$ and: if $j<j^*$ and $l<\omega$ then $\alpha^j_l=
\beta_m \Rightarrow \n(j)= \n\ \&\ p^j_l=\zeta$ ($=h^t_0(\n(j))(l))\in
Q_{\beta_m}$). 

In this case $\dom(q_{m+1})=\dom(q_m)\cup\{\beta_m\}$, and 
$$
q_{m+1}(\beta)=\left\{\begin{array}{ll} q_m(\beta) & \quad\mbox{\it if
}\beta<\beta_m (\mbox{so }\beta\in \dom(q_m))\\
\zeta & \quad\mbox{\it if } \beta=\beta_m
\end{array}
\right.
$$
\medskip

\noindent {\em Case C:} $m+1$ and for some $\n< \n^t$, we have $\n\in
\dom(h^t_1)$ and: $\alpha^j_l=\beta_m \Rightarrow \n(j)=\n$.

Work first in $\V[G_{P_{\beta_m}}]$, $q_m\in G_{P_{\beta_m}}$,
$G_{P_{\beta_m}}$ generic over $\V$. The sets
$$
\{\lim(\name{p}^j_l[G_{P_{\beta_m}}]): \alpha^j_l=\beta_m\mbox{ (and
}l\in [n^t_{k}, n^t_{k+1}), j<j^*)\}
$$
are subsets of $(^\omega2)^{[h^t_2(\n)]}=\{\eta\in {}^\omega2:
h^t_2(\n)\vartriangleleft \eta\}$, we can define an equivalence relation
$E_m$ on $(^\omega2)^{[h^t_2(\n)]}$:
$$
\nu_1 E_m \nu_2 \quad\mbox{\bf iff }\quad\nu_1\in\lim(\name{p}^j_l[
G_{P_{\beta_m}}]) \equiv \nu_2\in\lim(\name{p}^j_l[G_{P_{\beta_m}}])
$$
$$
\mbox{ whenever }\alpha^j_l=\beta_m.
$$
Clearly $E_m$ has finitely many equivalence classes, call them $\langle
Z^m_i:i< i^*_m\rangle$, all are Borel (sets of reals) hence they are
measurable; w.l.o.g.\ $\Leb(Z^m_i)=0\Leftrightarrow i\in [i^\otimes_m,
i^*_m)$, so clearly $i^\otimes_m>0$. For each $i<i^\otimes_m$ there is
$r=r_{m, i}\in \name{Q}_{\beta_m}[G_{P_{\beta_m}}]$ such that
$$
\lim(p^j_l[G_{P_{\beta_m}}])\supseteq Z^m_i \Rightarrow r\geq
p^j_l[G_{P_{\beta_m}}],
$$
$$
\lim(p^i_l[G_{P_{\beta_m}}])\cap Z^m_i=\emptyset\ \Rightarrow\ (\lim r)\cap
(\lim p^i_l[G_{P_{\beta_m}}])=\emptyset.
$$
We can also find a rational $a_{m, i}\in (0, 1)_{\bbr}$ such that
$$
a_{m, i} < \Leb(Z^m_i)/ 2^{\lgx(h^t_2(\n))} < a_{m, i} +\varepsilon/(2i^*_m).
$$
We can find $q'_m\in G_{P_{\beta_m}}$, $q_m\leq q'_m$ such that $q'_m$
forces all this information (so for $\name{Z}^m_i$, $\name{r}_{m, i}$ we
shall have $P_{\beta_m}$--names, but $a_{m, i}$, $i^\otimes_m$, $i^*_m$ are
actual objects). We then can find rationals $b_{m, i}\in (a_{m, i}, a_{m,
i}+\varepsilon/2)$ such that $\sum\limits_{i< i^\otimes_m} b_{m, i}=1$. Now
we throw a dice choosing $i_m< i^\otimes_m$ with the probability of $i_m=i$
being $b_{m, i}$ and define $q_{m+1}$ as:
$$
\dom(q_{m+1})=\dom(q'_m)\cup \{\beta_m\}
$$ 
$$
q_{m+1}(\beta)=\left\{\begin{array}{ll}
q'_m(\beta) &\quad \mbox{\it if }\beta<\beta_m \mbox{ (so }\beta\in
\dom(q'_m))\\
\name{r}_{m, i_m} & \quad \mbox{\it if }\beta=\beta_m
\end{array}\right.
$$
An important point is that this covers all cases (and in Case $B$ the choice
of $(j, l)$ is immaterial) as for each $\beta_m$ there is a unique $\n<
\n^t$ and $l$ such that $\eta_{\beta_m}\restriction w^t = \eta^t_{\n,
l}$ (see Definition \ref{K3} clause (b) and Definition \ref{T} clause
(i)). Basic probability computations (for $n^t_{k+1}-n^t_k$ independent 
experiments) show that for each $j$ coming from clause (a), by the law of
large numbers the probability of successes is $> 1-1/j^*$, successes meaning
$q_{m^*}\Vdash ``k\in \name{C}_j$'' (remember if $j$ comes from clause (b)
we always succeed).
\end{proof}

\begin{remark}
In the definition of $t\in\cT$ (i.e.\ \ref{T}) we can add $\eta^t_{\n,
\omega}\in {}^{w^t}2$ (i.e.\ replace $\langle\eta^t_{\n,l}:l<\omega\rangle$
by $\langle\eta^t_{\n,l}:l\leq\omega\rangle$) and demand 
\begin{enumerate}
\item[(l)] if $\zeta\in w^t$ then for every $n<\omega$ large enough,
$\zeta\in \eta^t_{\n,\omega}\equiv\zeta\in\eta^t_{\n,\omega}$,
\end{enumerate}
and in Definition \ref{K3} clause (d) use $\bar{\alpha}=\langle\alpha_l:
l<\omega\rangle$ {\em but} this does not help here.
\end{remark}

\begin{lemma}
\label{2.11}
1) Assume 
\begin{enumerate}
\item[(a)] $\bar{Q}\in {\mathcal K}^{3}$, $\bar{Q}=\langle
P_{\alpha},\name{Q}_{\beta},A_{\beta}, \mu_\beta,\name{r}_{\beta},
\eta_{\beta},(\name{\Xi}_{\alpha}^{t})_{t\in T}:\alpha \leq \alpha^{*},
\beta <\alpha^{*}\rangle$,
\item[(b)] $A\subseteq \alpha^{*}$, $\kappa <|A|<\lambda$,
\item[(c)] $\eta \in (\tka)^{V}\setminus\{\eta_{\beta}:\beta <\alpha^{*}\}$,
\item[(d)] $(\forall\alpha\in A)[|A_\alpha|<\kappa\ \Rightarrow\
A_\alpha\subseteq A]$ and $P'_{A}\lesdot P_{\alpha^{*}}$, $\name{Q}=
\name{Q}^{A,\bar Q}$ is the $P_{\alpha^*}$-name from \ref{2.1}(F)$(\beta)$
and  
$$
\mbox{if }\ \ t\in \cT\ \ \mbox{ then }\ \ \name{\Xi}^t_{\alpha^*}\restriction
\V^{P_A}\mbox{ is a $P_A$-name,}
$$
\end{enumerate}
{\em Then} there is 
$$
\bar{Q}^{+}=
\langle P_{\alpha},\name{Q}_{\beta},A_{\beta}, \mu_\beta, \name{\tau}_{\beta},
\eta_{\beta},(\name{\Xi}_{\alpha}^{t})_{t\in T}:\alpha \leq \alpha^{*}+1,
\beta <\alpha^{*}+1\rangle
$$ 
from $\cK^3$, extending $\bar{Q}$ such that 
$\name{Q}_{\alpha^{*}}=\name{Q}$, $A_{\alpha^{*}}=A$,
$\eta_{\alpha^{*}}=\eta$.

\noindent 2) If clauses (a)+(b)+(c) of part one hold {\em then} we can find
$A'$ such that: $A\subseteq A' \subseteq \alpha^*$, $|A'|\leq
(|A|+\kappa)^{\aleph_0}$ (which is $<\lambda$ by Hypothesis \ref{2.0}) and
such that $\bar Q$, $A'$, $\eta$ satisfy (a)+(b)+(c)+(d).
\end{lemma}

\begin{proof}
1) As before the problem is to define $\name{\Xi}^{t}_{\alpha^{*}+1}$.  We
have to satisfy clause (g) of Definition \ref{K3} for each fixed $t\in
\cT$. Let $\n^*$ be the unique $\n<\n^t$ such that $\eta\restriction
w^t=\eta^t_{\n, l}$. If $\n^*\in \dom(h^t_0)$ or $\langle \eta^t_{\n^*, l}:
l<\omega\rangle$ not constant or there is no such $\n^*$ then we have
nothing to do. So assume that $\alpha_{l}=\alpha^*$ for $l<\omega$,
$\eta^{t}_{\n^*,l}=\eta\restriction w^t$ (for $l<\omega$). Let $\Gamma$ be
the set of all pairs $( \name{r},\langle \name{r}_{l}: l<\omega \rangle)$
which satisfy the assumption $(**)^{\bar{Q}}_{\name{r},\langle\name{r}_l:l
<\omega\rangle}$ of \ref{K3} clause (g). In $\V^{P_{\alpha^{*}+1}}$ we have
to choose $\name{\Xi}^{t}_{\alpha^{*}+1}$ taking care of all these
obligations.  We work in $\V^{P_{\alpha^{*}}}$. By assumption (d) and Claim
\ref{1.6} it suffices to prove it for $\V^{P_A}$ so $Q_{\alpha^*}$ is 
$\Random^{\V^{P_A}}$ (see \ref{notl}(7)). By \ref{1.7} it is enough to prove
condition (B) of \ref{1.7}. Suppose it fails. Then there are $\langle B_m:
m<m(*)\rangle$ a partition of $\omega$ from $\V^{P_A}$, for simplicity
$\Xi^t_{\alpha^*}(B_m)>0$ for $m<m(*)$, and $(\name{r}^{i},\langle
\name{r}^{i}_{l}:l<\omega\rangle)\in \Gamma$ and $\n(i)=\n^*<\n^{t}$ for
$i<i^{*}<\omega$ and $\varepsilon^*>0$ and $r\in Q_{\alpha^{*}}$ which
forces the failure (of (B) of \ref{1.7}) for these parameters; (the
$\varepsilon^*$ comes from \ref{1.7}). W.l.o.g.\ $r$ forces that
$\name{r}^{i}\in \name{G}_{Q_{\alpha^{*}}}$ for $i<i^{*}$ (otherwise we can
ignore such $\name{r}^{i}$ as nothing is demanded on them in (g) of
\ref{K3}). So $r\geq r^{i}$ for $i<i^{*}$. 

By the assumption, for each $i<i^{*}$ we have: for each $r' \geq r$ (hence
$r'\geq r^i$ and $r' \in \Random$) and $i<i^{*}$ we have:
$$
Av_{\Xi^t_{\alpha^*}}(\langle a^i_k(r'): k<\omega\rangle)\geq 
(1-h^{t}_{1}(\n^*))
$$
where (see \ref{K3}(g)($\boxtimes$)) we let
$$
a^i_k(r')=a_k(\name{r},\langle\name{r}_l:l<\omega\rangle,\bar{n}^t)=
\frac{1}{n^t_{k+1}-n^t_{k}}\sum\limits_{l\in [n^t_k, n^t_{k+1})}
\frac{\Leb(\lim(r') \cap \lim(r^{i}_{l}))}{\Leb(\lim(r'))}.
$$ 
By \ref{1.7} it suffices to prove the following

\begin{lemma}
\label{2.15new}
Assume $\Xi$ is a finitely additive measure, $\langle B_0,\ldots,B_{m^*-1}
\rangle$ a partition of $\omega$, $\Xi(B_m)=a_m$, $i^*<\omega$ and $r$,
$r^i_l\in \Random$ for $i<i^*$, $l<\omega$ and $\bar{n}^*=\langle n^*_i:
i<\omega\rangle$, $n^*_i<n^*_{i+1}<\omega$ are such that 
\begin{enumerate}
\item[$(*)$] for every $r'\in \Random$, $r'\geq r$ and $i<i^*$ we have
$Av_\Xi(\langle a^i_k(r'): k<\omega\rangle)\geq b_i$ where
$$
a^i_k(r')=a^i_k(r',\langle r^i_l:l<\omega\rangle,\bar{n}^*) =\frac{1}{n^*_{k
+1}- n^*_k}\sum\limits^{n^*_{k+1}-1}_{l=n^*_k}\frac{\Leb(\lim(r')\cap
\lim(r^i_l))}{\Leb(\lim(r'))}.
$$
{\em Then} for each $\varepsilon>0$, $k^*<\omega$ there is a finite
$u\subseteq \omega\setminus k^*$ and $r'\geq r$ such that:
\begin{enumerate}
\item[(1)] $a_m-\varepsilon < |u\cap B_m|/|u| < a_m+\varepsilon$, for
$m<m^*$
\item[(2)] for each $i<i^*$ we have
$$
\frac{1}{|u|}\sum\limits_{k\in u} \frac{|\{l: n^*_k\leq l<n^*_{k+1}\mbox{
and }r'\geq r^i_l\}|}{n^*_{k+1}- n^*_k} 
$$
is $\geq b_i-\varepsilon$.
\end{enumerate}
\end{enumerate}
\end{lemma}

\begin{proof}
Let for $i < i^*$, $m< m^*$ and $r'\geq r$ (from $\Random$):
$$
c_{i, m}(r')= Av_{\Xi\restriction B_m}(\langle a^i_k(r'):
k\in B_m\rangle)\in [0, 1]_{\bbr}.
$$
So clearly 
\begin{enumerate}
\item[$(*)_1$] for $r'\geq r$ (in $\Random$)
$$
\begin{array}{ll}
b_i\leq Av_\Xi (\langle a^i_k(r'): k<\omega\rangle) &=\sum\limits_{m<m^*}
Av_{\Xi\rest B_m}(\langle a^i_k(r'): k\in B_m\rangle)\ \Xi(B_m)\\
\ & =\sum\limits_{m<m^*} c_{i, m}(r')
a_m.
\end{array}
$$
\end{enumerate}
There are $r^*\geq r$ and a sequence $\bar{c}=\langle c_{i,m}: i<i^*, m<
m^*\rangle$ such that: 
\begin{enumerate}
\item[$(*)_2$ (a)] $c_{i,m}\in [0, 1]_{\bbr}$,
\item[(b)] $\sum\limits_{m<m^*} c_{m,i} a_m \geq b_i$,
\item[(c)] for every $r'\geq r^*$ there is $r^{\prime\prime} \geq r^\prime$
such that
$$
(\forall i< i^*) (\forall m< m^*) [ c_{i, m}- \varepsilon < c_{i,
m}(r^{\prime\prime})< c_{i, m}+\varepsilon].
$$
\end{enumerate}
[Why? Let $k^*<\omega$ be such that $1/k^*<\varepsilon/(10\cdot l^*\cdot
m^*)$ (so $k^*>0$). Let $\Gamma=\{\bar{c}:\bar{c}=\langle c_{i,m}: i<i^*,\
m<m^*\rangle,\ c_{i,m}\in [0,1]_{\mathbb R}$ and $k^*c_{i,m}$ is an integer
and $\sum\limits_{m<m^*} c_{i,m}a_m>b_i\}$. Clearly $\Gamma$ is finite and
let us list it as $\langle\bar{c}^s:s<s^*\rangle$. We try to choose by
induction on $s\leq s^*$ a condition $r_s\in\Random$ such that $r_0=r$,
$r_s\leq r_{s+1}$, and for no $r''\geq r_{s+1}$ do we have 
\[(\forall i<i^*)(\forall m<m^*)[c^s_{i,m}-\varepsilon<c_{i,m}(r'')<
c^s_{i,m}+\varepsilon].\]
For $s=0$ we have no problem. If we succeed to arrive to $r_{s^*}$, for
$i<i^*$, $m<m^*$ we can define $c^*_{i,m}\in \{l/k^*:l\in\{0,\ldots,k^*\}\}$
such that $c_{i,m}(r_{s^*})\leq c^*_{i,m}<c_{i,m}(r_{s^*})+\varepsilon/(10
\cdot l^*\cdot m^*)$. By $(*)_1$ we have $b_i\leq\sum\limits_{m<m^*} c_{i,m}
(r^*_{s^*})a_m$. Clearly 
\[\sum_{m<m^*} c^s_{i,m} a_m\geq \sum_{m<m^*} c_{i,m}(r^*_{s^*}) a_m\]
so $\bar{c}^*=\langle c^*_{i,m}: i<i^*,\ m<m^*\rangle\in\Gamma$, hence for
some $s<s^*$, $\bar{c}^*=\bar{c}^s$. But then $r^*$ contradicts the choice
of $r_{s+1}$. Also by the above $\Gamma\neq\emptyset$. So we necessarily are
stuck at some $s<s^*$, i.e.\ cannot find $r_{s+1}$ as required. This means
that $r_s,\bar{c}^s$ as needed in $(*)_2$, so $r^*,\bar{c}$ as required
exist.]

Let $k^*<\omega$ be given. Now choose $s^*<\omega$ large enough and try to
choose by induction on $s\leq s^*$, a condition $r_s\in\Random$ and
natural numbers $(m_s, k_s)$ (flipping coins along the way) such that:
\begin{enumerate}
\item[$(*)_3$ (a)] $r_0=r^*$,
\item[(b)] $r_{s+1}\geq r_s$,
\item[(c)] $c_{i, m}-\varepsilon < c_{i, m}(r_s) < c_{i, m}+\varepsilon$ for
$i<i^*$, $m< m^*$,
\item[(d)] $k_s> k^*$, $k_{s+1} > k_s$,
\item[(e)] $k_s\in B_{m_s}$.
\end{enumerate}
In stage $s$, given $r_s$, we define $r_{s+1}$, $i_s$, $m_s$, $k_s$ as
follows: choose $m_s<m^*$ randomly with the probability of $m_s=m$ being
$a_m$. Next we can find a finite set $u_s\subseteq B_{m_s}$ with $\min(u_s)>
\max\{k^*+1, k_{s_1}+1: s_1<s\}$ such that
\begin{enumerate}
\item[$(*)$] if $i<i^*$ then $c_{i,m_s}-\varepsilon/2<\frac{1}{|u_s|}
\sum\limits_{k\in u_s} a^i_k(r_s) < c_{i, m_s}+\varepsilon/2$
\end{enumerate}
We define an equivalence relation $\e_s$ on $\lim(r_s)$:
$$
\eta_1\e_s\eta_2 \mbox{ iff }(\forall i<i^*)(\forall k\in u_s)
(\forall \ell\in[n^*_k, n^*_{k+1}))[\eta_1\in \lim(r^i_\ell)\equiv
\eta_2\in \lim(r^i_\ell)].
$$ 
The number of $\e_s$--equivalence classes is finite, and if $Y\in \lim(r_s)/
\e_s$ satisfies $\Leb(Y)>0$ choose $r_{s, Y}\in \Random$ such that
$\lim(r_{s, Y})\subseteq Y$ and $r_{s,Y}$ satisfies clause (c) of $(*)_3$
(possible by clause (c) of $(*)_2$). Now choose $r_{s+1}$ among $\{r_{s, Y}:
Y\in \lim(r_s)/\e_s \mbox{ and }\Leb(Y)>0\}$, with the probability of $r_{s,
Y}$ being $\Leb(Y)/\Leb(\lim(r_s))$. Lastly choose $k_s\in u_s$, with all
$k\in u_s$ having the same probability.

Now the expected value, assuming $m_s=m$, of
$$
\frac{1}{n^*_{k_s+1} - n^*_{k_s}} \times |\{\ell: n^*_{k_s}\leq \ell <
n^*_{k_s+1}\mbox{ and }r_{s+1}\geq r^i_\ell\}|
$$ 
belongs to the interval $(c_{i, m}- \varepsilon/2, c_{i, m}+\varepsilon/2)$,
because of the expected value of 
$$
\frac{1}{|u_s|} \sum\limits_{k\in u_s} \frac{1}{n^*_{k+1} -
n^*_k}\times |\{\ell: n^*_k\leq \ell < n^*_{k+1}\mbox{ and } r_{s+1}\geq
r^i_\ell\}|
$$
is in this interval (as 
\[\sum\{\Leb(Y):Y\in\lim(r_s)/\e_s,\ r_{s,Y}\geq r^i_l\}=\frac{\Leb(
\lim(r_s)\cap\lim(r^i_l))}{\Leb(\lim(r_s))}\]
and see the choice of $a^i_k(-)$).

Let $r'=r_{s^*}$, $u=\{k_s: s\leq s^*\}$. Hence the expected value of
$$
\frac{1}{|u|} \sum\limits_{k\in u} \frac{1}{n^*_{k+1} - n^*_k}\times |\{
\ell: n^*_k \leq \ell < n^*_{k+1}\mbox{ and }r'\geq r^i_\ell\}|
$$
is $\geq \sum a_m (c_{i, m} - \varepsilon/2)\geq b_i -\varepsilon/2$.

As $s^*$ is large enough with high probability (though just positive
probability suffices), $(r_{s^*}, \{k_s: s<s^*\})$ are as required for
$(r', u)$; note: we do not know the variance but we have a bound for
it not depending on $s$.

2) Straightforward.\end{proof}\end{proof}

The following is needed later to show that there are enough cases of the
Definition of $t$ with clause (g) of Definition \ref{K3} being non trivial
(i.e. $(**)$ there holds).

\begin{lemma}
\label{2.15}
Assume
\begin{enumerate}
\item[(a)] $\Xi$ is a finitely additive measure on $\omega$ and $b\in
(0,1]_{\bbR}$,
\item[(b)] $n^*_k<\omega$ (for $k<\omega$), $n^*_k<n^t_{k+1}$, and
$\lim(n^*_{k+1}- n^t_k)=\infty$,
\item[(c)] $r^*$, $r_l\in \Random$ are such that: 
\begin{enumerate}
\item[$(*)$] $(\forall l<\omega)[\Leb(\lim(r^*)\cap \lim(r_l))/
\Leb(\lim(r^*))\geq b].$
\end{enumerate}
\end{enumerate}
{\em Then} for some $r^\otimes\geq r^*$ we have:
\begin{enumerate}
\item[$\otimes_{r^\otimes}$] for every $r'\geq r^\otimes$ we have
$Av_\Xi(\langle a(r', k): k<\omega\rangle)\geq b$ where:

$a_k(r')=a(r', k)=a_k(\lim r')$ and for $X\subseteq {}^\omega 2$ we
let
$$
a_k(X)=
\frac{1}{n^*_{k+1}-n^*_k}\sum\limits_{l\in [n^*_k,n^*_{k+1})}\frac{\Leb(X
\cap (\lim r_l))}{\Leb(X)}. 
$$
\end{enumerate}
\end{lemma}

\begin{proof} Let 
$$
\cI=\{r\in \Random:\ r\geq r^*, \mbox{ and }Av_\Xi(\langle a_k(r'): k<\omega
\rangle) <b\}.
$$
If $\cI$ is not dense above $r^*$ there is $r^\otimes \geq r^*$ (in
$\Random$) such that for every $r\geq r^\otimes$, $r\notin \cI$ so
$r^\otimes$ is as required, so assume toward contradiction that $\cI$ is not
dense above $r^*$. There is a maximal antichain $\cI_1=\{s_i:i<i^*\}
\subseteq\cI$ (maximal among those $\subseteq \cI$), now $\cI_1$ is a
maximal antichain above $r^*$ as $r\in \cI \Rightarrow r\geq r^*$ and the
previous sentence. Hence $\Leb(\lim r^*)= \sum\limits_{i<i^*}\Leb(\lim
s_i)$; of course $|i^*|\leq \aleph_0$ as $\Random$ satisfies c.c.c. so 
w.l.o.g.\ $i^*\leq \omega$.

For any $j< i^*$ let $s^j= \bigcup\limits_{i<j} s_i$, note then
$\lim(\bigcup\limits_{m<i} s_m) = \bigcup\limits_{m<i}\lim(s_m)$ and
$$
a_k(s^j)= a_k (\bigcup\limits_{m<i} s_m)= \sum\limits_{i<j}
\frac{\Leb(s_i)}{\Leb(\bigcup\limits_{m<j} s_m)} a_k(s_i)
$$
hence
$$
\begin{array}{ll}
Av_\Xi(\langle a_k(s^j): k<\omega\rangle) &= Av_\Xi(\langle
a_k(\bigcup\limits_{m<j} s_m): k<\omega\rangle) \\
\ & = \sum\limits_{i<j}
\frac{\Leb(s_i)}{\Leb(\bigcup\limits_{m<j} s_m)}\times Av_\Xi(\langle
a_k(s_i): k<\omega\rangle) \\
\ & \leq \frac{\Leb(\lim (s_0))}{\Leb(\lim
(\bigcup_{i<j}s_i))}(b-\varepsilon) + \sum\limits_{0<i<j}
\frac{\Leb(\lim(s_i))}{\Leb(\lim(\bigcup_{m<j}s_m))} b\\
\ & = b- \Leb(\lim(s_0))\cdot \varepsilon
\end{array}
$$
where $\varepsilon= b - Av_{\Xi}( \langle a_k(s_0): k< \omega\rangle)$
so $\varepsilon >0$.

Let $j$ be large enough such that $\frac{\Leb(\lim(r^*)\setminus\lim(s^j)
)}{\Leb(\lim(r^*))}< \Leb(\lim(s_0))\cdot \varepsilon$. So
$$
\begin{array}{l}
Av_\Xi(\langle a_k(r^*): k<\omega\rangle)= \\
\qquad\frac{\Leb(\lim(r^*)\setminus \lim(s^j))}{\Leb(\lim(r^*))}
Av_\Xi(\langle a_k(\lim(r^*)\setminus \lim(s^j)): k<\omega\rangle) \\
\qquad\quad + \frac{\Leb(\lim(s^j))}{\Leb(\lim(r^*))} Av(\langle a_k(s^j):
k<\omega\rangle)\\
\qquad \leq \frac{\Leb(\lim(r^*)\setminus \lim(s^j))}{\Leb(r^*)} \times 1 +
\frac{\Leb(\lim(s^j))}{\Leb(\lim(r^*))} \times (b -
\Leb(\lim(s_0))\varepsilon)\\
\qquad < \Leb(\lim(s_0))\cdot \varepsilon + (b- \Leb(\lim(s_0))\cdot
\varepsilon)= b
\end{array}
$$
contradicting assumption (c). \end{proof}

\begin{claim}
\label{2.14}
Assume 
\begin{enumerate}
\item[(a)] 
$\bar{Q}\in {\mathcal K}^{3}$, $\bar{Q}=\langle P_{\alpha},\name{Q}_{\beta},
A_{\beta}, \mu_\beta,\name{r}_{\beta},\eta_{\beta},(\name{\Xi}_{\alpha}^{
t})_{t\in T}:\alpha \leq \alpha^{*},\beta <\alpha^{*}\rangle$, 
\item[(b)] $A\subseteq \alpha^*$ and $|A|<\kappa$ and $\mu<\kappa$ are such
that $\beta\in A\ \&\ |A_\beta|<\kappa\ \Rightarrow\ A_\beta\subseteq A$,
\item[(c)] $\eta\in {}^\kappa 2\setminus \{\eta_\beta: \beta<\alpha^*\}$,
\item[(d)] $\name{Q}$ is a $P_{\alpha^*}$-name of a forcing notion with set
of elements $\mu$, and is really definable in $\V[\langle \name{\tau}_\alpha:
\alpha\in A\rangle]$ from $\langle\name{\tau}_\alpha: \alpha\in A\rangle$
and parameters from $\V$.
\end{enumerate}
{\em Then} there is 
$$
\bar Q^+=\langle P\alpha, \name{Q}_\alpha, A_\beta, \mu_\beta,
\name{\tau}_\beta, \eta_\beta, (\name{\Xi}_\alpha^t)_{t\in \cT}: \alpha\leq
\alpha^*+1, \beta<\alpha^*+1\rangle
$$
from $\cK^3$ extending $\bar Q$ such that $\name{Q}_{\alpha^*}=\name{Q}$,
$A_{\alpha^*}=A$, $\eta_{\alpha^*}=\eta$ and $mu_{\alpha^*}=\mu$.
\end{claim}

\begin{proof} Straight. \end{proof}

\begin{remark}
If $Q$ is the Cohen forcing we can make one step toward $\{A\subseteq\omega:
\Xi^t_{\alpha^*+1}(A)=1\}$ being a selective filter but not needed at
present. 
\end{remark}

\section{Continuation of the proof of Theorem \protect\ref{A}}

We need the following lemma.

\begin{lemma}
\label{epsilon}
Suppose that $\bar{\varepsilon}=\langle \varepsilon_{l}:l<\omega\rangle$ is
a sequence of positive reals and $\bar Q \in {\mathcal K}^{3}$ has length
$\alpha$. The following set ${\mathcal I}_{\bar{\varepsilon}}\subseteq
P_{\alpha}$ is dense:
\[\begin{array}{ll}
{\mathcal I}_{\bar{\varepsilon}} = \{p\in P'_{\alpha}: &\mbox{there are }
m \mbox{ and } \alpha_{l},\nu_{l}\;(\mbox{for }l<m) \mbox{ such that }\\
\ & (a)\ \dom(p)=\{\alpha_{0}, \ldots ,\alpha_{m-1}\},\
\alpha_{0}>\alpha_{1}>\ldots >\alpha_{m-1},\\
\ & (b) \mbox{ if }|Q_{\alpha_{l}}|<\kappa,\mbox{ then }
p(\alpha_{l}) \mbox{ is an ordinal },\\ 
\ & (c) \mbox{ if } Q_{\alpha_{l}}\mbox{ is a partial random, then }
\forces_{P_{\alpha_{l}}}``p(\alpha_{l})\subseteq(^\omega2)^{[\nu_{l}]}\\ 
\ &\quad\mbox{ and }\Leb(\lim(p(\alpha_{l})))\geq
(1-\varepsilon_{l})/2^{\lgx(\nu_{l})}\mbox{\rm ''}\}\\
\end{array} \]
\end{lemma}

\begin{proof} By induction on $\alpha$ for all possible $\bar{\varepsilon}$.
\end{proof}

\begin{discussion}
\label{3.1A}
1) By the previous sections it follows that it is enough to prove that if
$\bar{Q} \in {\mathcal K}^{3}$, $P_{\alpha}=\Lim(\bar{Q})$, {\em then} in
$\V^{P_{\alpha}}$ the following sufficient condition holds:
\begin{enumerate}
\item[$(**)_{\bar{Q}}$] In $\V^{P_{\alpha}}$:\quad there is no perfect tree
$T \subseteq\fs$, $m\in \omega$ and $E\in [\lambda]^{\kappa^{+}}$ such that
$T \subseteq\tree_{m}[\bar{a}^{\alpha}]$ for all $\alpha \in E$.
\end{enumerate}

\noindent 2) Note that if we just want to prove $\forces_{P_\alpha}\mbox{``}
{\Gb}\leq \kappa\mbox{''}$ life is easier: $\Xi^t_\alpha$ is a zero-one
measure (so essentially an ultrafilter) and we interpret for
$\alpha<\lambda$, the forcing notion $Q_\alpha$ as $(^{\omega>}\omega,
\vartriangleleft)$ with generic real $\name{\eta}_\alpha$ and replace below
$(**)_{\bar Q}$ by
\begin{enumerate}
\item[$(**)_{\bar Q}^+$] in $V^{P_\alpha}$ there is no $\eta^*\in
{}^\omega\omega$ such that $\{\alpha<\lambda: (\forall \ell<\omega)
(\name{\eta}_\alpha(\ell)\leq \eta^*(\ell)\}$ has cardinality $\geq
\kappa^+$.
\end{enumerate}
In the proof below, $\name{T}$ is replaced by $\name{\eta}$, and
$p'_\zeta(\alpha_\zeta)$ is ${s^*}\conc \langle \zeta\rangle$.

\noindent 3) We can make the requirements on the $\Delta$-system stronger:
make it indiscernible also over some $A\subseteq \alpha$ of cardinality
$<\kappa$, where $\name{T}$ is a $P_A$-name, $p^*\in P_A$, and w.l.o.g.  the
heart is $\subseteq A$.

\noindent 4) Here the existence of $h^t_2$ help; we can use \ref{epsilon}
with $\sum\limits_{\ell<\omega}\varepsilon_\ell$ very small.
\end{discussion}

\begin{lemma}
\label{3.2}
If $P_{\alpha}=\Lim(\bar Q)$, $\alpha=\lgx(\bar Q)$ and 
$\bar Q \in {\mathcal K}^{3}$, {\em then} $(**)_{\bar Q}$ from
\ref{sufficient}. 
\end{lemma} 

\begin{proof}  
Suppose that $p^{*} \forces_{P_\alpha} `` \name{T},m,\name{E} \mbox{ \rm
form a counterexample to } (**)_{\bar{Q}}\mbox{\rm ''}$, wlog $p^*\in
P'_\alpha$. Let $\bar{\varepsilon}= \langle \varepsilon_{l}:l<\omega\rangle$
be such that $\varepsilon_{l}\in (0,1)_{\mathbb R}$ and $\sum\limits_{l<
\omega}\sqrt{\varepsilon_{l}} <1/10$. For each $\zeta<\kappa^{+}$ let
$p_{\zeta}\geq p^{*}$ be such that $p_{\zeta}\in {\mathcal I}_{\bar{
\varepsilon}}$ ($\subseteq P'_\alpha$) witnessed by $\langle
\nu^{\zeta}_\beta: \beta\in\dom (p_\zeta)\mbox{ and }|Q_\beta|\geq\kappa
\rangle$ (on $\cI_{\bar{\varepsilon}}$ see \ref{epsilon}) and 
$$
p_{\zeta} \forces_{P_\alpha} {\rm ``}\alpha_{\zeta} \mbox{ \rm is the }
\zeta\mbox{\rm -th element of }\name{E}\mbox{\rm ''}.
$$
So clearly $\alpha_\zeta< \lambda$. W.l.o.g., by thinning out, we can assume
that: 
\begin{itemize}
\item $\dom(p_{\zeta})=\{\gamma^{\zeta}_{i}:i<i^{*}\}$ with $\gamma^\zeta_i$
increasing with $i$, let $v^{\zeta}_{0}=\{i<i^{*}:|Q_{\gamma^{\zeta}_{i}}|<
\kappa\}$, then $v^{\zeta}_{0}=v_{0}$ is fixed for all $\zeta<\kappa^+$,
and let $v_{1}=i^{*}\setminus v_{0}$, 
\item $\dom(p_{\zeta})$ $(\zeta <\kappa^{+})$ form a $\Delta$--system, with 
the heart $\Delta$, so $\Delta\supseteq \dom(p^*)$,
\item $\alpha_{\zeta} \in \dom(p_{\zeta})$, $\alpha_{\zeta}=\gamma^{\zeta}_{z}$
for a fixed $z<i^{*}$,
\item $(\dom(p_{\zeta}),\Delta, <)$ are isomorphic for $\zeta<\kappa^{+}$,
\item if $i\in v_{0}$, then $p_{\zeta}(\gamma^{\zeta}_{i})=\gamma_{i}$,
for $\zeta<\kappa^{+}$,
\item if $i\in v_{1}$, then $\nu^{\zeta}_{\gamma^\zeta_i}=\nu_{i}$, 
(recall $\nu^{\zeta}_{\gamma^\zeta_i}\in \fs$ is given by the definition
of ${\mathcal I}_{\bar{\varepsilon}}$),
\item $p_{\zeta}(\alpha_{\zeta})=s^{*}$ for $\zeta<\kappa^{+}$ with $s^{*}=
\langle(n_{l},a_{l}):l<m^{*}\rangle$, w.l.o.g. $m^*>m$ (where $m$ is from
``the counterexample to $(**)_{\bar{Q}}$'') and $m^*>10$, 
\item for each $i<i^*$ the sequence $\langle \gamma^\zeta_i:\zeta<\kappa^+
\rangle$ is constant or strictly increasing,
\item the sequence $\langle \alpha_\zeta:\zeta< \kappa^+\rangle$ is
with no repetitions (as if $p_{\zeta_1}$, $p_{\zeta_2}$ are compatible
and $\zeta_1<\zeta_2<\lambda$ then $\alpha_{\zeta_1}\neq\alpha_{\zeta_2}$). 
\end{itemize}
Now we are interested only in the first $\omega$ conditions, i.e., we
consider $\zeta < \omega$.  For every such $\zeta$ let $p_{\zeta}' \geq
p_{\zeta}$ be such that $\dom(p_{\zeta}')=\dom(p_{\zeta})$,
$p_{\zeta}'(\gamma)=p_{\zeta}(\gamma)$ except for $\gamma=\alpha_{\zeta}$ in
which case we extend $p_{\zeta}(\alpha_{\zeta})=s^{*}$ in the following
way. We put $\lgx(p_{\zeta}'(\alpha_{\zeta}))=\lgx(s^{*})+1=m^{*}+1$,
$p_{\zeta}'(\alpha_{\zeta})={s^{*}}\conc \langle (j_{\zeta}^0,a_{\zeta})
\rangle$. Before we define $j^0_\zeta$, $a_{\zeta}$ choose an increasing
sequence of integers $\bar{s}=\langle s_{l}:l<\omega\rangle$, $s_{0}=0$,
such that $s_{k+1}-s_{k}=|[{}^{j_{k}}2]^{2^{j_{k}}(1-8^{-m^{*}})}|$ (i.e.\
it is the number of subsets of $^{j_k}2$ with $2^{j_k}(1-8^{-m^*})$
elements), where $j^{*}=3n_{m^{*}-1}+1$ (i.e. we define $j^*$ from the first
coordinate in the last pair in $s^*$) and we let $j_k=j^*+k!!$, and let
$j^0_\zeta=j_k$ when $\zeta\in [s_k, s_{k+1})$.  Now for $\zeta \in
[s_{k},s_{k+1})$ define $a_{\zeta}$ such that
\[\{a_{\zeta}:\zeta\in [s_{k},s_{k+1})\}=[{}^{j_k}2]^{2^{j_{k}}(1-
8^{-m^{*}})}\]
(so necessarily without repetitions). For $\varepsilon^{*}>0$ we define a
$P_\alpha$-name by 
$$
\name{A}_{\varepsilon^{*}}=\{k<\omega:|\{\zeta\in [s_{k},s_{k+1}):
p'_{\zeta}\in \name{G}_{P_{\alpha}}\}|/(s_{k+1}-s_{k})>\varepsilon^{*}\}.
$$
For the proof of \ref{3.2} we need: 
\begin{subclaim}
\label{3.3.1}
There is a condition $p^{\otimes}\geq p^{*}$ which forces that for some
$\varepsilon^{*}>0$ the set $\name{A}_{\varepsilon^{*}}$ is infinite .
\end{subclaim}

\begin{proof} 
Choose $\varepsilon^*>0$ small enough. First we define a suitable blueprint
$t\in \cT$,  
$$
t=(w^t, \n^{t}, \m^t, \bar \eta^t, h^{t}_{0},h^{t}_{1},h^{t}_{2}, \bar n^t ).
$$
Let 
$$
\begin{array}{r}
w^{t}=\{\min\{\beta<\kappa:\eta_{\gamma^{\zeta(1)}_{i(1)}}(\beta)\neq
\eta_{\gamma^{\zeta(2)}_{i(2)}}(\beta)\}:\zeta(1), \zeta(2)<\omega
\mbox{ and }\ \\
i(1), i(2)<i^*\mbox{ and } \gamma^{\zeta(1)}_{i(1)}\neq \gamma^{\zeta(2)
}_{i(2)}\ \}. 
\end{array}
$$
Let $\n^{t}=i^{*}$, $\dom(h^t_{0})=v_{0}$, $\dom(h^{t}_{1})=\dom(h^t_{2})=
v_{1}$ and $n^t_l=s_l$. If $\n\in v_{0}$, then $h^t_{0}(\n)(l)=\gamma_{\n}$
and $\eta^{t}_{\n,\zeta}=\eta_{\gamma^{\zeta}_{\n}}\restriction w^t$. If
$\n\in v_{1}$, then $h^t_1(\n)=\varepsilon_{n}$, $h^t_2(\n)=\nu_{\n}$.

We now define a condition $p^\otimes$, it will be in $P_\alpha$,
$\dom(p^\otimes)=\Delta$, $p^*\leq p^\otimes$; remember $\dom(p^*)\subseteq
\Delta$ as for each $\zeta$ we have $p^*\leq p_\zeta$.  If $\gamma\in
\Delta$ then for some $\n<\n^t$, we have $\bigwedge\limits_{\zeta <\omega}
\gamma^\zeta_{\n}=\gamma$.  If $\n\in v_0$ we let
$p^\otimes(\gamma)=h^t_0(\n)$, so trivially in $\V^{P_\gamma}$
$$
p^\otimes(\gamma)\Vdash_{Q_\gamma}{\rm `` }\name{\Xi}^t_{\gamma+1}(\{\zeta<
\omega: h^t_0(\n)\in G_{Q_\gamma}\})=1\mbox{ if }\n\in \dom(h^t_0) (=v_0)
\mbox{ ''}
$$

If $\n\in v_1$, then define a $P_{\gamma}$-name for a member of $Q_{\gamma}$
as follows. Consider $\name{r}^\n_{\zeta}=\name{p}_{\zeta}'(\gamma)$ for
$\zeta<\omega$.  Let $\name{r}$ be the member $(^\omega 2)^{[h^t_2(n)]}$ of
$Q_\gamma$. Working in $\V^{P'_{A_\alpha}}$, by Lemma \ref{2.15} there is
$\name{r}^{*}_{\gamma}\geq \name{r}$ from $\name{Q}_\gamma$ such that for
every $r' \geq r^{*}_\gamma$ in $Q_{\gamma}$ we have 
\begin{enumerate}
\item[$(**)_{r', \varepsilon}$] 
$Av_{\Xi^t_{\alpha}}(\langle a^\n_k(r'): k<\omega\rangle)\geq (1-h^t_1(\n))=
(1-\varepsilon_{\n})$
where
$$
a^\n_k(r') =:
\frac{1}{n^t_{k+1}-n^t_k}\sum\limits_{l\in [n^t_k,n^t_{k+1})}\frac{\Leb(
\lim(r')\cap\lim(r^\n_l))}{\Leb(\lim(r'))}.
$$
\end{enumerate}
Hence the assumption of condition $(g)$ in Definition \ref{K3} holds, hence
in $\V^{P_{\gamma}}$ we have:
$$
\begin{array}{r}
r^*_\gamma\forces_{Q_\gamma}\mbox{\rm `` }
Av_{\Xi^{t}_{\gamma+1}}(\langle |\{\ell\in [n^t_{k+1}- n^t_k):p_{\ell}(
\gamma)\in\name{G}_{Q_{\gamma}}\}|/ (n^t_{k+1}-n^t_k):k\in \omega\rangle)\\
\geq 1-\varepsilon_n\mbox{ ''}.
\end{array}
$$
So there is a $P_\gamma$-name $\name{r}^*_\gamma$ of such a condition.
In this case let $p^\otimes(\gamma)=\name{r}^*_\gamma$, so we have finished
defining $p^\otimes$, clearly it has the right domain. 

Now suppose that $\n<\n^t$, $\n\in v_1$ is such that $\gamma^{\zeta}_{\n}
\not\in\Delta$. Define $\bar{\beta}=\langle \beta_{\zeta}:\zeta
<\omega\rangle$, $\beta_{\zeta}=\gamma^{\zeta}_{\n}$. Then $\bar{\beta}$
satisfies $(t,\n)$ for $\bar{Q}$. By our assumption the assumption of clause
$(f)$ in Definition \ref{K3} is satisfied, hence in $\V^{P_{\alpha}}$, for
any $\varepsilon>0$:
$$
\forces_{P_\alpha}\mbox{\rm `` }\name{\Xi}^t_{\alpha}\Big(\{k\!: \frac{|\{l\in
[n^t_k, n^t_{k+1})\!: p_l(\gamma^l_n)\in \name{G}_{Q_{\gamma^l_n}}\}|}{(n^t_{k+
1}- n^t_k)}\geq (1-\varepsilon_n)\cdot (1-\varepsilon)\}\Big)=1\mbox{ ''}. 
$$
Now for each $\n\in v_1$, as 
$$
(1-\varepsilon_\n)\cdot (1-\varepsilon)\leq Av_{\Xi^t_\alpha}(\langle\frac{|
\{\ell: n^t_k\leq \ell<n^t_{k+1}\mbox{ and }r^\n_\ell\in G_{P_\alpha}\}|}{
n^t_{k+1}- n^t_k}:k<\omega\rangle)
$$
and $\varepsilon>0$ was arbitrary, clearly 
\begin{enumerate}
\item[$(*)_\n$] in $\V^{P_\alpha}$, $\sqrt{2\varepsilon_n}\geq
\Xi^t_\alpha\Big(\{ k<\omega: 1-\sqrt{2\varepsilon_n}\geq \frac{|\{\ell:
n^t_k\leq \ell < n^t_{k+1}\mbox{ and }r^\n_\ell\in
G_{P_\alpha}\}|}{n^t_{k+1}- n^t_k}\}\Big)$.
\end{enumerate}
Let
$$\name{A}^\prime_{\varepsilon^*}= \{k<\omega:\mbox{ if }\zeta\in
[s_k, s_{k+1})\mbox{ and } i\in v_o\mbox{ then }p_\zeta\restriction
\{\gamma^\zeta_i\}\in \name{G}_{P_\alpha}\},
$$
clearly $\Xi^t_\alpha(\name{A}^\prime_{\varepsilon^*})=1$.

Let $\varepsilon^*< 1- \sum\limits_{\n} \sqrt{2\varepsilon_\n}$ and
$\varepsilon^*>0$, so
$$
\begin{array}{r}
\name{A}_{\varepsilon^*}\cup(\omega\setminus A^\prime_{\varepsilon^*})
=\{k<\omega:\varepsilon^*<\frac{|\bigcap\limits_{\n\in v_1}\{l\in [n^t_k,
n^t_{k+1})\mbox{ and } p'_\zeta(\gamma^\zeta_{\n})\in
G_{Q_{\gamma_{\n}}}\}|}{(n^t_{k+1}-n^t_k)}\}\\
\\
\supseteq \{ k<\omega: \mbox{ if }\n\in v_1\mbox{ then }\frac{|\{\ell:
n^t_k\leq \ell< n^t_{k+1}\mbox{ and }r^{\n}_\ell\in G_{P_\alpha}\}|}{
n^t_{k+1}- n^t_k}\geq 1-\sqrt{2\varepsilon_\n}\}\\
\\
= \omega\setminus \bigcup\limits_{\n\in v_1}\{k<\omega:
\frac{|\{\ell: n^t_k\leq \ell< n^t_{k+1} \mbox{ and }r^\n_\ell\in
G_{P_\alpha}\}|}{n^t_{k+1}- n^t_k}< 1- \sqrt{2\varepsilon_\n}\},
\end{array}
$$
hence\ \ $\Xi^t_\alpha(\name{A}_{\varepsilon^*}\cup (\omega\setminus
\name{A}^\prime_{\varepsilon^*}))\geq 1-\sum\limits_{\n\in 
v}\sqrt{2\varepsilon_\n} \geq \varepsilon^*>0$, but
$$
\name{\Xi}^t_\alpha(\omega\setminus \name{A}^\prime_{\varepsilon^*})=1-
\name{\Xi}^t_\alpha(\name{A}^\prime_{\varepsilon^*})=1-1=0
$$
hence necessarily $A_{\varepsilon^*}$ is infinite.

This suffice for \ref{3.3.1}. \end{proof}

Let $p^{\otimes}$ be as in the claim \ref{3.3.1}, let $G_{P_{\alpha}} $ be a
generic subset of $P_\alpha$ to which $p^\otimes$ belongs and we shall work
in $\V[G_{P_\alpha}]$. So $A=\name{A}_{\varepsilon^{*}}[G]$ is infinite. For
$k\in A$, let $b_{k}=\{\zeta\in [s_k,s_{k+1}):p_{\zeta}'\in G_{P_\alpha}\}$.
We know that $|b_{k}|> (s_{k+1}-s_{k})\times\varepsilon^{*}$. Note that if
$k\in A$, then $T\cap {}^{j_{k}}2 \subseteq \bigcap_{\zeta\in b_{k}}
a_{\zeta}$ as $\lgx (s^*)=m^*>m$. To reach a contradiction it is enough to
show that for infinitely many $k\in A$ there is a bound on the size of
$c_{k}=T \cap {}^{j_{k}}2$ which does not depend on $k$.

Now $|b_k|/(s_{k+1}-s_k)$ is at most the probability that if we choose a
subset of $^{(j_k)}2$ with $2^{j_k}(1-8^{-m^*})$ elements, it will include
$T\cap {}^{(j_k)}2$; now if $k\in A$ this probability has a lower bound
$\varepsilon^*$ not depending on $k$, and this implies a bound on $|T\cap
{}^{(j_k)}2|$ not depending on $k$. More formally, for a fixed $k<\omega$ we
have:
\[\begin{array}{lcl}
|b_{k}| & = & |\{a_{\zeta}:\zeta\in [s_k,s_{k+1}),\: \zeta\in
b_{k}\}| \\ 
& \leq & |\{a_{\zeta}:\zeta\in [s_k,s_{k+1}),\: (T\cap {}^{j_{k}}2)
\subseteq a_{\zeta}\}| \\
& \leq & |\{a\subseteq (^{j_{k}}2): (T\cap {}^{j_{k}}2)
\subseteq a \mbox{ \rm and } |a|=2^{j_{k}}(1-8^{-m^{*}})\}| \\
& = & |\{a\subseteq (^{j_{k}}2)\setminus (T\cap {}^{j_{k}}2):
|a|=2^{j_{k}}\times 8^{-m^{*}}\}| 
\end{array}\]

$$
= \left(
\begin{array}{l}
2^{j_k}- |T\cap {}^{(j_k)}2|\\
2^{j_k}\cdot 8^{-m^*}
\end{array}
\right)
$$
similarly $s_{k+1}- s_k= \left(\begin{array}{l} 2^{j_k}\\
2^{j_k}\cdot 8^{-m^*}
\end{array}
\right).$ Hence
$$\begin{array}{ll}
|b_k|/(s_{k+1}- s_k)& \leq \left(\begin{array}{l}
2^{j_k}-|T\cap {}^{j_k}2|\\ 2^{j_k}\cdot
8^{-m^*}\end{array}\right)\Big/ \left(\begin{array}{l} 2^{j_k}\\
2^{j_k}\cdot 8^{-m^*}\end{array}\right) \\
\ & = \prod\limits_{i<|T\cap {}^{(j_k)}2|} (2^{j_k} - 2^{j_k}
8^{-m^*}-i)/ \prod\limits_{i< |T\cap {}^{(j_k)}2|}(2^{j_k} - i)\\
\ & = \prod\limits_{i< |T\cap {}^{(j_k)}2|}\left( 1-\frac{2^{j_k}
8^{-m^*}}{2^{j_k}-i}\right).
\end{array}
$$
Let $i_k(*) =\min\{|T\cap {}^{(j_k)}2|, 2^{j_k -1}\}$ so 
$$
\begin{array}{ll}
\varepsilon^*\leq \frac{|b_k|}{s_{k+1}- s_k} & \leq
\prod\limits_{i<|T\cap {}^{j_k}2|} (1-
\frac{2^{j_k}8^{-m^*}}{2^{j_k}-i})\\
\ & \leq \prod\limits_{i< i_k(*)} (1-
\frac{2^{j_k} 8^{-m^*}}{2^{j_k}})= (1- 8^{-m^*})^{i_k(*)}.
\end{array}
$$
So we can find a bound to $i_k(*)$ not depending on $k$:
$$
i_k(*)\leq (\log(1/\varepsilon^*)/ \log(1/(1-8^{-m^*}),
$$ 
remember $m^*> 10$
so $1-8^{-m^*}\in (0, 1)_{\bbr}$. So for $k$ large enough,
$$
|T\cap{}^{(j_k)}2| = i_k(*)\leq
\log(1/\varepsilon^*)/\log(1/(1-8^{-m^*})).
$$ 
This finishes the proof. \end{proof}

\begin{theorem}
\label{3.2new}
Under Hypothesis \ref{2.0} there is $\bar Q\in \cK^3$, $\lgx(\bar
Q)=\chi=\delta^*$ (if clause $(\alpha)$ of \ref{2.0}(b) holds) or $\lgx(\bar
Q)=\chi\times\chi\times\lambda^+$ (if clause $(\beta)$ of \ref{2.0}(b)
holds) such that in $\V^{P_{\lim \bar Q}}$ we have $MA_{<\kappa}+\cov(\nul)
=\lambda$. 
\end{theorem}

\begin{proof} 
First assume clause $(\alpha)$ of \ref{2.0}. By \ref{notl}(2) and
\ref{notl}(6) it suffices to find an iteration 
$$
\langle P_\alpha, \name{Q}_\beta, A_\beta, \mu_\beta, \name{\tau}_\beta,
\eta_\beta, (\name{\Xi}^t_\alpha)_{t\in \cT}: \alpha\leq \chi,
\beta<\chi\rangle\in\cK^3
$$ 
(see definition \ref{K3}) satisfying clauses (a)+(b)+(c) of
\ref{notl}(2)+(6) (as the only property missing, $\cov(\nul)\leq \lambda$, 
holds by \ref{sufficient} + \ref{3.2}.

Let $\cK_3^-=\{\bar Q\in \cK_3: \lgx(\bar Q)<\chi\}$. 

Now choose $\bar Q^\xi\in \cK_3^-$ for $\xi<\chi$ increasing with $\xi$ (see
definition \ref{2.9}) by induction on $\chi$. Now if $\cf(\xi)>\aleph_0$ use
\ref{2.9A}(3), if $\cf(\xi)=\aleph_0$ use \ref{2.10}. Bookkeeping give us
sometimes a case $\name{Q}$ of \ref{notl}(6)(c) as assignment, we can find
suitable $A\subseteq \lgx({\bar Q}^\xi)$ by \ref{notl}(4) and then apply
\ref{2.14} to get ${\bar Q}^{\xi+1}$. For other $\xi$, bookkeeping gives us
a case of \ref{notl}(2)(b) as assignment $A\subseteq \lgx({\bar Q}^\xi)$,
such that $|A|<\lambda$. Now we apply \ref{2.11}(2) (with $\bar Q$, $A$
there standing for ${\bar Q}'$, $A$ here) and get $A'$ as there. Now apply
\ref{2.11}(1) with ${\bar Q}'$, $A'$ here standing for $\bar Q$, $A'$ here
standing for $\bar Q$, $A$ there (and $\eta$ any member of
$^{\kappa}2\setminus \{\eta_\beta: \beta<\lgx({\bar Q}')\}$ possible as
$\lgx({\bar Q}')<\chi$ as $\bar Q'\in \cK^3$) and get $\bar Q^{\xi+1}$
(corresponding to $\bar Q^+$ there).

Second assume clause $(\beta)$ of \ref{2.0}(b). We just should be more
careful in our bookkeeping, particularly in the beginning let
$\langle\eta_\alpha:\alpha<\chi\times\chi\times\lambda^+\rangle$ be an
enumeration of ${}^\kappa 2$ with no repetition.
\end{proof}

\bibliographystyle{lit-plain}
\bibliography{listb,lista,listf,listx}
\shlhetal

\end{document}